\documentclass[11pt]{article}
\usepackage{amsmath, upgreek}
\usepackage{amssymb}
\usepackage{amscd}
\usepackage{xspace}
\usepackage{verbatim}

\newcommand{\mysection}[1]{
\section{#1}\setcounter{equation}{0}}
\title{\bf Local and global properties of solutions of quasilinear Hamilton-Jacobi equations}
\author{{\bf Marie-Fran\c{c}oise Bidaut-V\'eron\footnote{\noindent Laboratoire de Math\'{e}matiques et Physique Th\'{e}orique, CNRS UMR 7350,
Facult\'{e} des Sciences, 37200 Tours France. E-mail: veronmf@univ-tours.fr},} \\{\bf Marta Garcia-Huidobro \footnote{\noindent 
Departamento de Matematicas, Pontifica Universidad Catolica de Chile
Casilla 307, Correo 2, Santiago de Chile. E-mail: mgarcia@mat.puc.cl}}\\
 {\bf Laurent V\'eron \footnote{\noindent 
Laboratoire de Math\'{e}matiques et Physique Th\'{e}orique, CNRS UMR 7350,
Facult\'{e} des Sciences, 37200 Tours France. E-mail: veronl@univ-tours.fr}}\\[2mm]
}

\date{}
\begin{document}
 \maketitle


\newcommand{\txt}[1]{\;\text{ #1 }\;}
\newcommand{\tbf}{\textbf}
\newcommand{\tit}{\textit}
\newcommand{\tsc}{\textsc}
\newcommand{\trm}{\textrm}
\newcommand{\mbf}{\mathbf}
\newcommand{\mrm}{\mathrm}
\newcommand{\bsym}{\boldsymbol}
\newcommand{\scs}{\scriptstyle}
\newcommand{\sss}{\scriptscriptstyle}
\newcommand{\txts}{\textstyle}
\newcommand{\dsps}{\displaystyle}
\newcommand{\fnz}{\footnotesize}
\newcommand{\scz}{\scriptsize}
\newcommand{\be}{\begin{equation}}
\newcommand{\bel}[1]{\begin{equation}\label{#1}}
\newcommand{\ee}{\end{equation}}
\newcommand{\eqnl}[2]{\begin{equation}\label{#1}{#2}\end{equation}}
\newcommand{\barr}{\begin{eqnarray}}
\newcommand{\earr}{\end{eqnarray}}
\newcommand{\bars}{\begin{eqnarray*}}
\newcommand{\ears}{\end{eqnarray*}}
\newcommand{\nnu}{\nonumber \\}
\newtheorem{subn}{\name}
\renewcommand{\thesubn}{}
\newcommand{\bsn}[1]{\def\name{#1}\begin{subn}}
\newcommand{\esn}{\end{subn}}
\newtheorem{sub}{\name}[section]
\newcommand{\dn}[1]{\def\name{#1}}   
\newcommand{\bs}{\begin{sub}}
\newcommand{\es}{\end{sub}}
\newcommand{\bsl}[1]{\begin{sub}\label{#1}}
\newcommand{\bth}[1]{\def\name{Theorem}
\begin{sub}\label{t:#1}}
\newcommand{\blemma}[1]{\def\name{Lemma}
\begin{sub}\label{l:#1}}
\newcommand{\bcor}[1]{\def\name{Corollary}
\begin{sub}\label{c:#1}}
\newcommand{\bdef}[1]{\def\name{Definition}
\begin{sub}\label{d:#1}}
\newcommand{\bprop}[1]{\def\name{Proposition}
\begin{sub}\label{p:#1}}
\newcommand{\R}{\eqref}
\newcommand{\rth}[1]{Theorem~\ref{t:#1}}
\newcommand{\rlemma}[1]{Lemma~\ref{l:#1}}
\newcommand{\rcor}[1]{Corollary~\ref{c:#1}}
\newcommand{\rdef}[1]{Definition~\ref{d:#1}}
\newcommand{\rprop}[1]{Proposition~\ref{p:#1}}
\newcommand{\BA}{\begin{array}}
\newcommand{\EA}{\end{array}}
\newcommand{\BAN}{\renewcommand{\arraystretch}{1.2}
\setlength{\arraycolsep}{2pt}\begin{array}}
\newcommand{\BAV}[2]{\renewcommand{\arraystretch}{#1}
\setlength{\arraycolsep}{#2}\begin{array}}
\newcommand{\BSA}{\begin{subarray}}
\newcommand{\ESA}{\end{subarray}}
\newcommand{\BAL}{\begin{aligned}}
\newcommand{\EAL}{\end{aligned}}
\newcommand{\BALG}{\begin{alignat}}
\newcommand{\EALG}{\end{alignat}}
\newcommand{\BALGN}{\begin{alignat*}}
\newcommand{\EALGN}{\end{alignat*}}
\newcommand{\note}[1]{\textit{#1.}\hspace{2mm}}
\newcommand{\Proof}{\note{Proof}}
\newcommand{\qeda}{\hspace{10mm}\hfill $\square$}
\newcommand{\qed}{\\
${}$ \hfill $\square$}
\newcommand{\Remark}{\note{Remark}}
\newcommand{\modin}{$\,$\\[-4mm] \indent}
\newcommand{\forevery}{\quad \forall}
\newcommand{\set}[1]{\{#1\}}
\newcommand{\setdef}[2]{\{\,#1:\,#2\,\}}
\newcommand{\setm}[2]{\{\,#1\mid #2\,\}}
\newcommand{\mt}{\mapsto}
\newcommand{\lra}{\longrightarrow}
\newcommand{\lla}{\longleftarrow}
\newcommand{\llra}{\longleftrightarrow}
\newcommand{\Lra}{\Longrightarrow}
\newcommand{\Lla}{\Longleftarrow}
\newcommand{\Llra}{\Longleftrightarrow}
\newcommand{\warrow}{\rightharpoonup}
\newcommand{
\paran}[1]{\left (#1 \right )}
\newcommand{\sqbr}[1]{\left [#1 \right ]}
\newcommand{\curlybr}[1]{\left \{#1 \right \}}
\newcommand{\abs}[1]{\left |#1\right |}
\newcommand{\norm}[1]{\left \|#1\right \|}
\newcommand{
\paranb}[1]{\big (#1 \big )}
\newcommand{\lsqbrb}[1]{\big [#1 \big ]}
\newcommand{\lcurlybrb}[1]{\big \{#1 \big \}}
\newcommand{\absb}[1]{\big |#1\big |}
\newcommand{\normb}[1]{\big \|#1\big \|}
\newcommand{
\paranB}[1]{\Big (#1 \Big )}
\newcommand{\absB}[1]{\Big |#1\Big |}
\newcommand{\normB}[1]{\Big \|#1\Big \|}
\newcommand{\produal}[1]{\langle #1 \rangle}

\newcommand{\thkl}{\rule[-.5mm]{.3mm}{3mm}}
\newcommand{\thknorm}[1]{\thkl #1 \thkl\,}
\newcommand{\trinorm}[1]{|\!|\!| #1 |\!|\!|\,}
\newcommand{\bang}[1]{\langle #1 \rangle}
\def\angb<#1>{\langle #1 \rangle}
\newcommand{\vstrut}[1]{\rule{0mm}{#1}}
\newcommand{\rec}[1]{\frac{1}{#1}}
\newcommand{\opname}[1]{\mbox{\rm #1}\,}
\newcommand{\supp}{\opname{supp}}
\newcommand{\dist}{\opname{dist}}
\newcommand{\myfrac}[2]{{\displaystyle \frac{#1}{#2} }}
\newcommand{\myint}[2]{{\displaystyle \int_{#1}^{#2}}}
\newcommand{\mysum}[2]{{\displaystyle \sum_{#1}^{#2}}}
\newcommand {\dint}{{\displaystyle \myint\!\!\myint}}
\newcommand{\q}{\quad}
\newcommand{\qq}{\qquad}
\newcommand{\hsp}[1]{\hspace{#1mm}}
\newcommand{\vsp}[1]{\vspace{#1mm}}
\newcommand{\ity}{\infty}
\newcommand{\prt}{\partial}
\newcommand{\sms}{\setminus}
\newcommand{\ems}{\emptyset}
\newcommand{\ti}{\times}
\newcommand{\pr}{^\prime}
\newcommand{\ppr}{^{\prime\prime}}
\newcommand{\tl}{\tilde}
\newcommand{\sbs}{\subset}
\newcommand{\sbeq}{\subseteq}
\newcommand{\nind}{\noindent}
\newcommand{\ind}{\indent}
\newcommand{\ovl}{\overline}
\newcommand{\unl}{\underline}
\newcommand{\nin}{\not\in}
\newcommand{\pfrac}[2]{\genfrac{(}{)}{}{}{#1}{#2}}

\def\ga{\alpha}     \def\gb{\beta}       \def\gg{\gamma}
\def\gc{\chi}       \def\gd{\delta}      \def\ge{\epsilon}
\def\gth{\theta}                         \def\vge{\varepsilon}
\def\gf{\phi}       \def\vgf{\varphi}    \def\gh{\eta}
\def\gi{\iota}      \def\gk{\kappa}      \def\gl{\lambda}
\def\gm{\mu}        \def\gn{\nu}         \def\gp{\pi}
\def\vgp{\varpi}    \def\gr{\rho}        \def\vgr{\varrho}
\def\gs{\sigma}     \def\vgs{\varsigma}  \def\gt{\tau}
\def\gu{\upsilon}   \def\gv{\vartheta}   \def\gw{\omega}
\def\gx{\xi}        \def\gy{\psi}        \def\gz{\zeta}
\def\Gg{\Gamma}     \def\Gd{\Delta}      \def\Gf{\Phi}
\def\Gth{\Theta}
\def\Gl{\Lambda}    \def\Gs{\Sigma}      \def\Gp{\Pi}
\def\Gw{\Omega}     \def\Gx{\Xi}         \def\Gy{\Psi}

\def\CS{{\mathcal S}}   \def\CM{{\mathcal M}}   \def\CN{{\mathcal N}}
\def\CR{{\mathcal R}}   \def\CO{{\mathcal O}}   \def\CP{{\mathcal P}}
\def\CA{{\mathcal A}}   \def\CB{{\mathcal B}}   \def\CC{{\mathcal C}}
\def\CD{{\mathcal D}}   \def\CE{{\mathcal E}}   \def\CF{{\mathcal F}}
\def\CG{{\mathcal G}}   \def\CH{{\mathcal H}}   \def\CI{{\mathcal I}}
\def\CJ{{\mathcal J}}   \def\CK{{\mathcal K}}   \def\CL{{\mathcal L}}
\def\CT{{\mathcal T}}   \def\CU{{\mathcal U}}   \def\CV{{\mathcal V}}
\def\CZ{{\mathcal Z}}   \def\CX{{\mathcal X}}   \def\CY{{\mathcal Y}}
\def\CW{{\mathcal W}} \def\CQ{{\mathcal Q}}
\def\BBA {\mathbb A}   \def\BBb {\mathbb B}    \def\BBC {\mathbb C}
\def\BBD {\mathbb D}   \def\BBE {\mathbb E}    \def\BBF {\mathbb F}
\def\BBG {\mathbb G}   \def\BBH {\mathbb H}    \def\BBI {\mathbb I}
\def\BBJ {\mathbb J}   \def\BBK {\mathbb K}    \def\BBL {\mathbb L}
\def\BBM {\mathbb M}   \def\BBN {\mathbb N}    \def\BBO {\mathbb O}
\def\BBP {\mathbb P}   \def\BBR {\mathbb R}    \def\BBS {\mathbb S}
\def\BBT {\mathbb T}   \def\BBU {\mathbb U}    \def\BBV {\mathbb V}
\def\BBW {\mathbb W}   \def\BBX {\mathbb X}    \def\BBY {\mathbb Y}
\def\BBZ {\mathbb Z}

\def\GTA {\mathfrak A}   \def\GTB {\mathfrak B}    \def\GTC {\mathfrak C}
\def\GTD {\mathfrak D}   \def\GTE {\mathfrak E}    \def\GTF {\mathfrak F}
\def\GTG {\mathfrak G}   \def\GTH {\mathfrak H}    \def\GTI {\mathfrak I}
\def\GTJ {\mathfrak J}   \def\GTK {\mathfrak K}    \def\GTL {\mathfrak L}
\def\GTM {\mathfrak M}   \def\GTN {\mathfrak N}    \def\GTO {\mathfrak O}
\def\GTP {\mathfrak P}   \def\GTR {\mathfrak R}    \def\GTS {\mathfrak S}
\def\GTT {\mathfrak T}   \def\GTU {\mathfrak U}    \def\GTV {\mathfrak V}
\def\GTW {\mathfrak W}   \def\GTX {\mathfrak X}    \def\GTY {\mathfrak Y}
\def\GTZ {\mathfrak Z}   \def\GTQ {\mathfrak Q}

\font\Sym= msam10 
\def\SYM#1{\hbox{\Sym #1}}
\newcommand{\bdw}{\prt\Gw\xspace}

\date{}
\maketitle\medskip

\noindent{\small {\bf Abstract} We study some properties of the solutions of (E) $\;-\Gd_p u+|\nabla u|^q=0$ in a domain $\Gw \sbs \BBR^N$, mostly when $p\geq q>p-1$. We give a universal priori estimate of the gradient of the solutions with respect to the distance to the boundary. We give a full classification of the isolated singularities of the nonnegative solutions of (E), a partial classification of isolated singularities of the negative solutions. We prove a general removability result expressed in terms of some Bessel capacity of the removable set. We extend our estimates to equations on complete non compact manifolds satisfying a lower bound estimate on the Ricci curvature, and derive some Liouville type theorems.
}\smallskip

\noindent
{\it \footnotesize 2010 Mathematics Subject Classification}. {\scriptsize 35J92, 35J61, 35F21, 35B53
58J05}.\\
{\it \footnotesize Key words}. {\scriptsize p-Laplacian; a priori estimates; singularities; removable set; Bessel capacities; curvatures; convexity radius.
}
\vspace{1mm}
\hspace{.05in}

\tableofcontents
\medskip
\mysection{Introduction}
Let  $N\geq p>1$, $q>p-1$ and $\Gw \sbs \BBR^N$ ($N>1$) be a  domain. In this article we study some local and global properties of solutions of  
\be\label{A1}\left.\BA{ll}
-\Gd_p u+|\nabla u|^q=0
\EA\right.\ee
in $\Gw$, where $\Gd_pu:=\rm{div}\left(|\nabla u|^{p-2}|\nabla u|\right)$. The main questions we consider are the following:\medskip

 \noindent 1- A priori estimates and Liouville type theorems. \smallskip
 
 \noindent 2- Removability of singularities. \smallskip

 \noindent 3- Description of isolated singularities of solutions.\medskip
  
 Our technique allows us to handle both nonnegative and signed solutions. We will speak of a problem with {\it absorption} when we consider nonnegative solutions and a problem with {\it source} when we consider negative solutions (in which case we will often set $u=-\tilde u$). One of the main tools we use is a pointwise gradient estimate, valid for {\it any signed} solution of $(\ref{A1})$, \medskip
 
 \noindent {\bf Theorem A}. {\it Let $u\in C^1(\Gw)$ be a solution of $(\ref{A1})$ in $\Gw$. Then
 \bel{I-1}
 \abs{\nabla u(x)}\leq c_{N,p,q}(\dist(x,\prt\Gw))^{-\frac{1}{q+1-p}}
 \ee 
 for any $x\in\Gw$. If $\Gw=\BBR^N$, $u$ is a constant.}\medskip

In the case $p=2$ the existence of an upper bound of the gradient has first been obtained by Lasry and Lions \cite{LaLi} and then made explicit by Lions \cite{Lio}; the idea there was based upon the Bernstein technique. In \cite{NpV} Nguyen-Phuoc and   V\'eron rediscovered this upper bound by a slightly different method. Our method of proof is a combination of the Bernstein approach and the Keller-Osserman construction of radial supersolutions of the elliptic inequality satisfied by $\abs{\nabla u}^2$, a technique which will be fundamental for extension of this results in a geometric framework (see below).  \smallskip

 Concerning solutions of $(\ref{A1})$ in a domain $\Gw\subset\BBR^N$ we obtain that if $p\neq q$, any solution satisfies
   \bel{I-4}
\abs{u(x)}\leq c_{p,q}\left(\dist (x,\prt\Gw)^{\frac{p-q}{q+1-p}}-{\gd_*}^{\frac{p-q}{q+1-p}}\right)+
\max\{\abs{u(z)}:\dist (z,\prt\Gw)=\gd_*\}
 \ee 
if $\dist (x,\prt\Gw)\leq\gd_*$, where $\gd_*>0$ depends on the curvature of $\prt\Gw$;  when $p=q $  the formula holds provided the term $\dist (x,\prt\Gw)^{\frac{p-q}{q+1-p}}-{\gd_*}^{\frac{p-q}{q+1-p}}$ is replaced by $\ln (\dist (x,\prt\Gw)/\gd_*)$. 
In the case $p=2$ this estimate was a key element for the study of boundary singularity developed in \cite{NpV}. This aspect of equation $(\ref{A1})$ will be developed in a forthcoming article.
\medskip

In the study of singularities, we first give a general removability result concerning interior singularities. The general removability result given in \rth{remov}, is expressed in terms of the Bessel capacity $C_{1,\frac{q}{q+1-p}}$ relative to $\BBR^N$, and deals with locally renormalized solutions (see \rdef {renorm}).\medskip

\noindent {\bf Theorem B} {\it Let $p-1<q<p$ and $F\subset \Gw$ be a relatively closed set such that $C_{1,\frac{q}{q+1-p}}(F)=0$. Then any locally renormalized solution $u$ of $(\ref{A1})$ in $\Gw\setminus F$ can be extended as a locally renormalized solution in whole $\Gw$. When $u$ is nonnegative, $u$ is therefore a $C^1$ solution in whole $\Gw$. When $u$ is a signed renormalized solution, it is a $C^1$ solution provided $q<\frac{N(p-1)}{N-1}$.}\medskip

Our result is actually stronger and deals with locally renormalized solutions of 
\be\label{M1}\left.\BA{ll}
-\Gd_p u+|\nabla u|^q=\gm
\EA\right.\ee
in $\Gw\setminus F$ where $\gm$ is a measure which is absolutely continuous with respect to the $C_{1,\frac{q}{q+1-p}}$ capacity. \smallskip

A previous result in \cite{Phuc} shows the same result provided $u$ is a p-subharmonic function in $\Gw$, $\gm=0$ and $q>p$. An interesting counter part in \cite{Phuc0} asserts that if any weak solution 
of $(\ref{A1})$ in $\Gw\setminus F$ can be extended as a solution in $\Gw$, then $C_{1,\frac{q}{q+1-p}}(F)=0$. The construction therein is strongly associated with the solvability of problem $(\ref{M1})$ by the use the capacitary measure of $F$. When $p=2$, necessary and sufficient conditions for the existence of a solution of $(\ref{M1})$ in $\BBR^N$ can be found in  \cite{HMV}. \medskip

When $(\ref{A1})$ is replaced by 
\be\label{M2}\left.\BA{ll}
-\Gd_p u+\ge u^q=\gm
\EA\right.\ee
with $\ge=\pm1$ and $p>q>p-1$, deep existence results  of solutions $(\ref{M2})$ of have been obtained in  \cite{PhV1}, \cite{PhV2}, with $\ge=-1$ and \cite{B-V}, \cite{BV-H-V} with $\ge=1$. Their proofs (excepted for  \cite{B-V}) are strongly based upon fine study of Wolff potentials and their links with Bessel potentials. In the case $\ge=-1$, a necessary and sufficient condition for existence of a positive solution of $(\ref{M2})$ is obtained with an assumption of Lipschitz continuity of the measure $\gm\geq 0$ with respect to the $C_{p,\frac{q}{q+1-p}}$ capacity. In the case $\ge=1$, a sufficient condition for existence of a signed solution of $(\ref{M2})$  for a signed measure $\gm$ is obtained with the assumption that $\gm$ is absolutely continuous with respect to the $C_{p,\frac{q}{q+1-p}}$ capacity.

\medskip

If $K$ is reduced to a single point $0\in \Gw$, the threshold of removability of an isolated singularity corresponds to the exponent 
   \bel{I-5}
   q=q_c:=\frac{N(p-1)}{N-1}
\ee
 but the situation is different if we consider positive or negative solutions. \smallskip

If $q>p-1$, we set
\bel{I-7-}
\gb_q=\myfrac{p-q}{q+1-p}.
\ee
When $p-1<q<q_c$, there exists an explicit radial positive solution of $(\ref{A1})$ in $\BBR^{N}_*=\BBR^N\setminus\{0\}$
   \bel{I-6}
U(x)=U(|x|)=\gl_{N,p,q}|x|^{-\gb_q}
 \ee 
 where
    \bel{I-7}\BA{lll}
\gl_{N,p,q}=\gb^{-1}_q\left(\gb_q(p-1)+p-N\right)^{\frac{1}{q+1-p}}.
 \EA\ee 
When  $p=2$, Lions obtained in \cite{Lio} the description of isolated singularities of nonnegative solutions of $(\ref{A1})$ in the subcritical case $1<q<\frac{N}{N-1}$. We extend his result to the general case $1<p\leq N $ and provide a full classification of isolated singularities of nonnegative solutions :
 \medskip
 
  \noindent{\bf Theorem C }{\it Assume $p-1<q<q_c$ and $u\in C^1(\Gw\setminus\{0\})$ is a nonnegative solution of $(\ref{A1})$ in $\Gw\setminus \{0\}$. Then
 
 \noindent (i) either there exists $c\geq 0$ such that 
  \bel{I-8}
\lim_{x\to 0}\frac{u(x)}{\gm_p(x)}=c,
 \ee 
 where $\gm_p$ is the fundamental solution of the p-Laplacian defined by 
 $$\gm_p(x)=\frac{p-1}{N-p}\abs x^{\frac{p-N}{p-1}}\quad\text{ if }\;1<p<N\quad\text{ and }\;\gm_N(x)=-\ln\abs x.$$
  Furthermore $u$ satisfies
   \bel{I-8'}
-\Gd_pu+\abs{\nabla u}^q=c_{N,p}c\gd_0\qquad\text{in }\CD'(\Gw),
 \ee 
 
 \noindent (ii) or 
   \bel{I-9}
\lim_{x\to 0}\abs x^{\gb_q}u(x)=\gl_{N,p,q}
 \ee 
for some explicit positive constants $\gl_{N,p,q}$ and $\gb_q=\frac{p-q}{q+1-p}$.}\medskip

When $q\geq q_c$ the nonnegative solutions can be extended as $C^1$ functions. Concerning negative solutions there exists a radial negative singular solutions $V=-\tilde U$ of $(\ref{A1})$ in $\BBR^{N}_*$ with
   \bel{I-10}
\tilde U(x)=\tilde U(|x|)=\tilde \gl_{N,p,q}|x|^{-\gb_q},
 \ee 
 where
    \bel{I-11}
\tilde \gl_{N,p,q}=\gb^{-1}_q\left(N-p-\gb_q(p-1)\right)^{\frac{1}{q+1-p}}.
 \ee 
 In this case we obtain  a  partial classification of isolated singularities of negative solutions of $(\ref{A1})$ in $\Gw\setminus \{0\}$.\medskip
 
  \noindent{\bf Theorem D }{\it Assume $u$ is a negative $C^1$ solution of $(\ref{A1})$ in $\Gw\setminus \{0\}$. Then
 
  \noindent(i) {\it  When $p-1<q<q_c$ there exists $c\leq 0$ such that $(\ref{I-8})$ and $(\ref{I-8'})$ hold.}

 \noindent (ii) {\it  When $q> q_c$ $(\ref{I-8})$  and $(\ref{I-8'})$ hold with $c=0$. \smallskip
 
 \noindent Furthermore, when $u$ is  radial, there holds: \smallskip
 
 \noindent (iii)   When $q=q_c$, either
 
     \bel{I-12}
\lim_{x\to 0}\myfrac{\left(-\ln |x|\right)^{\frac{N-1}{p-1}}u(x)}{\gm_p(x)}=c_{N,p}<0
 \ee
 or $u$ is regular at $0$}. \smallskip
 
  \noindent (iv)   When $q>q_c$, 
 
     \bel{I-13}
\lim_{x\to 0}|x|^{\gb_q}u(x)=-\tilde \gl c_{N,p,q},
 \ee
 or $u$ is regular at $0$}.\medskip
 
 In the last section we obtain local and global estimates of solutions when $\BBR^N$ is replaced by a N-dimensional  Riemannian manifold $(M^{^{_N}}\!\!,g)$ and $-\Gd _p$ by the corresponding p-Laplacian $-\Gd_{g,p}$ in covariant derivatives. Our results  emphasize the role of the Ricci curvature of the manifold if $p=2$ and the sectional curvature if $p\neq 2$. In the case $1<p< 2$ we need to introduce the notion of {\it convexity radius} of a point $x\in M$, denoted by $r_M(x)$, which is supremum of the $r>0$ such that the geodesic ball $B_{r}(x)$ is strongly convex.\medskip

 \noindent {\bf Theorem E}. {\it Let $q>p-1>0$, $(M^{^{_N}}\!\!,g)$ be a  Riemannian manifold with Ricci curvature $Ricc_g$ and sectional curvature $Sec_g$ and $\Gw\subset M$ be a domain such that $Ricc_g\geq (1-N)B^2$ in $\Gw$. Assume  also $Sec_g\geq -\tilde B^2$ in $\Gw$ if $p>2$ or $r_M(z)\geq \dist (z,\prt\Gw)$ for any $z\in\Gw$ if $1<p<2$. Then any $C^1$ solution $u$ of 
\be\label{A1*}\left.\BA{ll}
-\Gd_{g,p}u+|\nabla u|^q=0
\EA\right.\ee
in $\Gw$ satisfies
%
 \bel{I-2}
 \abs{\nabla u(x)}^2\leq c_{N,p,q}\max\left\{B^{\frac{2}{q+1-p}},(1+B_pd(x,\prt\Gw))^{\frac{1}{q+1-p}}(d(x,\prt\Gw))^{-\frac{2}{q+1-p}}\right\},
 \ee 
 for any $x\in\Gw$, where $B_p=B+(p-2)_+\tilde B$ and $d(x,\prt\Gw)$ is now the geodesic distance of $x$ to $\prt\Gw$.}   \medskip
 
Notice that {\it $r_M(x)$ is always infinite when $Sec_g\leq 0$}.  Furthermore if for some $a\in M$ we have that  $r_M(a)=\infty$, then $r_M(x)=\infty$ for any $x\in M$; in this case we say that the convexity radius $r_M$ of $M$ is infinite. As a consequence we obtain \medskip
 
  \noindent {\bf Theorem F}. {\it Let $0<p-1<q$ and $(M^{^{_N}}\!\!,g)$ be a  complete noncompact Riemannian manifold such that $Ricc_g\geq (1-N)B^2$, ($B\geq 0$). Assume also $r_M=\infty$ if $1<p<2$ or
   \bel{I-2'}
\lim_{\dist(x,a)\to\infty}\myfrac{|Sec_g(x)|}{\dist(x,a)}=0
 \ee
 for some $a\in M$ if $p>2$. Then any solution $u$ of $(\ref{A1*})$ in $M$ satisfies
  \bel{I-2''}
 \abs{\nabla u(x)}\leq c_{N,p,q}B^{\frac{1}{q+1-p}}\qquad\forall x\in M.
 \ee }
 
  Since our estimate holds also in the case $p=q$,  we obtain \medskip
  
  \noindent {\bf Theorem G}. {\it   Assume $M$ satisfies the assumptions of Theorem F. Then any positive p-harmonic  function $v$ on $M$ satisfies
  \bel{I-3}
 v(a)e^{-\gk B\dist(x,a)}\leq v(x) \leq v(a)e^{\gk B\dist(x,a)}
 \ee 
 for any points $a,x$ in $M$, where $\gk=\gk(p,N)>0$}. \medskip 
 
 The case $p=2$, $B=0$ is due to Chen and  Yau (\cite{CY}). Kortschwar and Li \cite {KL} obtain a similar estimate but a with a global estimate of the sectional curvature which implies our assumption on the Ricci curvature. \medskip

 \noindent{\bf Aknowledgements} This article has been prepared with the support of the MathAmsud collaboration program 13MATH-02 QUESP. The first two authors were supported by Fondecyt grant N¡ 1110268. The authors are grateful to A. El Soufi for helpful discussions and to the referee for interesting references.

\mysection{A priori estimates in a domain of $\BBR^N$}

\subsection{The gradient estimates}
The next result is the extension to the  $p$-Laplacian of  a result obtained by Lions \cite{Lio} for the Laplacian. We denote by $d(x)$ the distance from $x\in\overline\Gw$ to $\prt\Gw$. 

\bprop{lets0} Assume $q>p-1$ and $u$ is a $C^1$ solution of $(\ref{A1})$ in a domain $\Gw$. Then
\bel{est0}
\abs{\nabla u(x)}\leq c_{N,p,q}(d(x))^{-\frac{1}{q+1-p}}\qquad\forall x\in \Gw.
\ee
\es
\Proof In any open subset $G$ of $\Gw$ where $|\nabla u|>0$ we write $(\ref{A1})$ under the form
\bel{grad1}
-\Gd u-(p-2)\myfrac{D^2u(\nabla u,\nabla u)}{\abs{\nabla u}^2}+\abs{\nabla u}^{q+2-p}=0
\ee
and we recall the formula
\bel{W1}\myfrac{1}{2}\Gd\abs{\nabla u}^2=\abs{D^2u}^2+\langle\nabla \Gd u,\nabla u\rangle. 
\ee
By Schwarz inequality
$$\abs{D^2u}^2\geq \myfrac{1}{N}(\Gd u)^2,
$$
hence we obtain
\bel{grad1'}
\myfrac{1}{2}\Gd\abs{\nabla u}^2\geq \myfrac{1}{N}(\Gd u)^2+\langle\nabla \Gd u,\nabla u\rangle.
\ee
Next, we write $z=\abs{\nabla u}^2$ and derive from $(\ref{grad1})$
\bel{grad2}
\Gd u=-\myfrac{(p-2)}{2}\myfrac{\langle \nabla z,\nabla u\rangle}{z}+z^{\frac{q+2-p}{2}},
\ee
thus
\bel{grad2'}\BA {l}
\langle\nabla\Gd u,\nabla u\rangle=-\myfrac{(p-2)}{2}\myfrac{\langle D^2z(\nabla u),\nabla u\rangle}{z}-\myfrac{(p-2)}{4}\myfrac{\abs{\nabla z}^2}{z}
\\[3mm]\phantom{------\langle\nabla\Gd u,\nabla u\rangle}
+\myfrac{(p-2)}{2}\myfrac{\langle \nabla z,\nabla u\rangle^2}{z^2}+
\myfrac{q+2-p}{2}z^{\frac{q-p}{2}}\langle\nabla z,\nabla u\rangle.
\EA\ee
From $(\ref{grad2})$ and $(\ref{grad2'})$
\bel{grad3}
(\Gd u)^2\geq \frac{1}{2}z^{q+2-p}-\myfrac{(p-2)^2}{4}\myfrac{\langle\nabla z,\nabla u\rangle^2}{z^2},
\ee
and from $(\ref{grad1'})$,
\bel{grad4}\BA {l}
\Gd z+(p-2)\myfrac{\langle D^2z(\nabla u),\nabla u\rangle}{z}\geq\myfrac{1}{N}z^{q+2-p}-\myfrac{(p-2)^2}{2N}\myfrac{\langle\nabla z,\nabla u\rangle^2}{z^2}
\\[3mm]\phantom{\Gd z--}
-\myfrac{(p-2)}{2}\myfrac{\abs{\nabla z}^2}{z}+(p-2)\myfrac{\langle \nabla z,\nabla u\rangle^2}{z^2}+
(q+2-p)z^{\frac{q-p}{2}}\langle\nabla z,\nabla u\rangle.
\EA\ee
Noticing that $\myfrac{\langle \nabla z,\nabla u\rangle^2}{z^2}\leq \myfrac{|\nabla z|^2}{z},$
and that for any $\ge>0$ 
$$z^{\frac{q-p}{2}}\abs{\langle\nabla z,\nabla u\rangle}\leq z^{\frac{q+2-p}{2}}\myfrac{\abs{\nabla z}}{\sqrt z}\leq \ge
z^{{q+2-p}}+\frac{1}{4\ge}\myfrac{\abs{\nabla z}^2}{z},
$$
we obtain that the right-hand side of $(\ref{grad4})$ is bounded from below by the quantity
$$\frac{1}{2N}z^{q+2-p}-D\myfrac{|\nabla z|^2}{z},
$$
where $D=D(p,q,N)>0$.
We define the operator
\bel{grad5}
v\mapsto \CA (v):=-\Gd v-(p-2)\myfrac{\langle D^2v(\nabla u),\nabla u\rangle}{\abs{\nabla u}^2}=-\sum_{i,j=1}^Na_{ij}v_{x_ix_j}
\ee
where the $a_{ij}$ depend on $\nabla u$;  since $\abs{\nabla u}^2=z$,  
\bel{grad6}\BA{ll}\displaystyle
\gth\abs\xi^2:=\min\{1,p-1\}\abs\xi^2\leq \sum_{i,j=1}^Na_{ij}\xi_i\xi_j\leq \max\{1,p-1\}\abs\xi^2:=\Gth\abs\xi^2,
\EA\ee
for all $\xi=(\xi_1,...,\xi_N)\in\BBR^N$. Consequently, $\CA$ is uniformly elliptic in $G$ and $z$ satisfies
\bel{grad7}
\CL(z):=\CA (z)+\frac{1}{2N}z^{q+2-p}- D\myfrac{\abs{\nabla z}^2}{z}\leq 0\qquad\text {in }\Gw.
\ee
Consider a ball $B_a(R)\subset\Gw$ and set $w(x)=\tilde w(\abs {x-a})=\gl(R^2-\abs {x-a}^2)^{-\frac{2}{q+1-p}}$. Put $r=\abs {x-a}$, then 
$$w_{x_i}=\myfrac{4\gl}{q+1-p}(R^2-\abs {x-a}^2)^{-\frac{2}{q+1-p}-1}x_i,$$
$$w_{x_ix_j}=\myfrac{4\gl}{q+1-p}(R^2-r^2)^{-\frac{2}{q+1-p}-1}\gd_{ij}+
\myfrac{8(3+q-p))\gl}{(q+1-p)^2}(R^2-r^2)^{-\frac{2}{q+1-p}-2}x_ix_j,
$$
therefore
$$\abs{\nabla w}^2=\myfrac{16\gl^2}{(q+1-p)^2}(R^2-r^2)^{-\frac{4}{q+1-p}-2}r^2,
$$

$$\myfrac{\abs{\nabla w}^2}{w}=\myfrac{16\gl}{(q+1-p)^2}(R^2-r^2)^{-\frac{2}{q+1-p}-2}r^2,
$$

$$w^{q+2-p}= \gl^{q+2-p}(R^2-r^2)^{-2\frac{q+2-p}{q+1-p}}=\gl^{q+2-p}(R^2-r^2)^{-\frac{2}{q+1-p}-2},
$$
and finally 
$$\CA(w)\geq -\myfrac{4\Gth\gl}{q+1-p}(R^2-r^2)^{-\frac{2}{q+1-p}-2}\left(NR^2+\left(\myfrac{3+q-p}{q+1-p}-N\right)r^2\right).
$$
At the end, using the fact that $r\leq R$, we obtain
$$\CL(w)\geq \myfrac{\gl}{2N}(R^2-r^2)^{-\frac{2}{q+1-p}-2}\left(\gl^{q+1-p}-cR^2\right),
$$
where $c=c_{N,p,q}$. We choose $\gl=(cR^2)^\frac{1}{q+1-p}$ and derive $\CL(w)\geq 0$. We take for $G$ a connected component of $\{x\in B_R(a):z(x)>w(x)\}$, thus $z(x)>0$ in $G$ and $\overline G\subset \overline B_R(a)$. If $x_0\in G$ is such that
$z(x_0)-w(x_0)=\max\{z(x)-w(x):x\in  G\}$, then $\nabla z(x_0)=\nabla w(x_0)$, $z(x_0)>w(x_0)>0$ and 
$$\CA(z(x_0)-w(x_0))+\frac{1}{2N}\left(z(x_0)^{q+2-p}-w(x_0)^{q+2-p}\right)-D\abs{\nabla z}^2\left(\myfrac{1}{z(x_0)}-\myfrac{1}{w(x_0)}\right)\leq 0,
$$
which contradicts the fact that all the terms are nonnegative with the exception of  $z(x_0)^{q+2-p}-w(x_0)^{q+2-p}$ which is positive. 
Therefore $z\leq w$ in $B_R(a)$. In particular 
\bel{grad8}
z(a)\leq w(a)=c'_{N,p,q}R^{-\frac{2}{q+1-p}}.
\ee
Letting $R\to d(x)$ yields $(\ref{est0})$.\qeda

\subsection{Applications}
The first estimate is a pointwise one for solutions with possible isolated singularities if $q\leq p$.

\bcor {ptws}Assume $q>p-1>0$, $\Gw$ is a domain containing $0$ and let $R^*>0$ be such that $d(0)\geq 2R^*$. Then for any $x\in B_{R^*}\setminus\{0\}$, and $0<R\leq R^*$, any  $u\in C^2(\Gw\setminus\{0\})$ solution of $(\ref{A1})$ in $\Gw\setminus\{0\})$ satisfies
\bel{estu0-1}
\abs{u(x)}\leq c_{N,p,q}\abs {|x|^{-\gb_q}-R^{-\gb_q}}+\max\{\abs{u(z)}:\abs z=R\},
\ee
if $p\neq q$, and 
\bel{estu0-2}
\abs{u(x)}\leq c_{N,p}\left(\ln \abs R-\ln \abs x\right)+\max\{\abs{u(z)}:\abs z=R\},
\ee
if $p=q$.
\es
\Proof Let $X=\frac{R}{\abs x}x$, then  $\dist(tx+(1-t)X,\prt B_R(X))=t\abs x+(1-t)R$ for any $0<t<1$, thus by $(\ref{est0})$ in $B_{R^*}\setminus\{0\}$
$$\BA {l}
\abs{u(x)}=\abs{u(X)+\myint{0}{1}\myfrac{d}{dt}u(tx+(1-t)X)dt}\\[4mm]\phantom{\abs{u(x)}}
\leq \abs{u(X)}+c_{N,p,q}\abs{x-X}\myint{0}{1}(t\abs x+(1-t)R)^{-\frac{1}{q+1-p}}dt.
\EA$$
By integration, we obtain $(\ref{estu0-1})$ or $(\ref{estu0-2})$. In the particular case where $p>q$ and $\abs x\leq \frac{R}{2}$, we obtain
\bel{estu0-3}
\abs{u(x)}\leq c_{N,p,q}|x|^{-\gb_q}+\max\{\abs{u(z)}:\abs z=R\}.
\ee
\qeda\medskip

The second estimate corresponds to solutions with an eventual boundary blow-up if $q\leq p$..

\bcor{unif}Assume $q>p-1>0$, $\Gw$ is a bounded domain with a $C^2$ boundary. Then there exists $\gd_*>0$ such that if we denote $\Gw_{\gd_*}:=\{z\in\Gw:d(z)\leq \gd_*\} $, any   $u\in C^2(\Gw) $ solution of $(\ref{A1})$ in $\Gw$ satisfies 
\bel{estu1}
\abs{u(x)}\leq c_{N,p,q}\abs {(d(x))^{-\gb_q}-\gd^{*{-\gb_q}}}+\max\{\abs{u(z)}:d(z)=\gd_*\}\quad\forall x\in \Gw_{\gd_*}
\ee
 if $p\neq q$ and 
\bel{estu2}
\abs{u(x)}\leq c_{N,p,q}\left(\ln \gd_*-\ln d(x)\right)+\max\{\abs{u(z)}:d(z)=\gd_*\}\quad\forall x\in \Gw_{\gd_*}
\ee
if $p=q$.
\es
\Proof We denote by $\gd_*$ the maximal $r>0$ such that any boundary point $a$ belongs to a ball $B_r(a_i)$ of radius $r$  such that $B_r(a_i)\subset \overline\Gw$ and to a ball $B_r(a_s)$ with radius $r$ too such that $B_r(a_s)\subset \overline\Gw^c$. If $x\in \Gw_{\gd_*}$, we denote by $\gs(x)$ its projection onto $\prt\Gw$ and by ${\bf n}_{\gs(x)}$ the outward normal unit  vector to $\prt\Gw$ at $\gs(x)$ and $z^*=\gs(x)-2\gd_*{\bf n}_{\gs(x)}$. Then 
$$u(x)=u(z^*)+\myint{0}{1}\myfrac{d}{dt}u(tx+(1-t)z^*)dt=\myint{0}{1}\langle\nabla u(tx+(1-t)z^*),x-z^*\rangle dt
$$
thus
$$\abs{u(x)}\leq \abs{u(z^*)}+c_{N,p,q}(\gd_*-d(x))\myint{0}{1}(td(x)+(1-t)d(z^*))^{-\frac{1}{q+1-p}}dt.
$$
Integrating this relation we obtain $(\ref{estu1})$ and $(\ref{estu2})$. \qeda\medskip

\noindent\Remark As a consequence of $(\ref{estu1})$ there holds for $p>q>p-1$
\bel{estu3}
u(x)\leq \left(c_{N,p,q}+K\max\{\abs{u(z)}:d(z)\geq\gd_*\}\right)(d(x))^{-\gb_q}\quad\forall x\in \Gw
\ee
where $K=({\rm diam} (\Gw))^{\frac{p-q}{q+1-p}}$, with the standard modification if $p=q$.\medskip

As a variant of \rcor{unif} we have estimate of solutions in an exterior domain
\bcor{unif2}Assume $q>p-1>0$, $R_0>0$ and let  $u\in C^2(B_{R_0}^c)$ be any solution of $(\ref{A1})$ in $B_{R_0}^c$. Then for any $R>R_0$ there holds
\bel{estu1-1}
\abs{u(x)}\leq c_{N,p,q}\abs{(\abs x-R_0)^{-\gb_q}-(R-R_0)^{-\gb_q}}+\max\{\abs{u(z)}:\abs z=R\}\quad\forall x\in B_{R}^c
\ee
 if $p\neq q$ and 
\bel{estu2-2}
\abs{u(x)}\leq c_{N,p,q}\left(\ln (\abs x-R_0)-\ln (R-R_0)\right)+\max\{\abs{u(z)}:\abs z=R\}\quad\forall x\in B_{R}^c
\ee
if $p=q$.
\es
\Proof The proof is a consequence of the identity
$$u(x)=u(z)+\myint{0}{1}\myfrac{d}{dt}u(tx+(1-t)z)dt=\myint{0}{1}\langle\nabla u(tx+(1-t)z),x-z\rangle dt
$$
where $z=\frac{R}{\abs x}x$. Since 
$$\abs{\nabla u(tx+(1-t)z)}\leq C_{N,p,q}(t\abs{x}+(1-t)R-R_0)^{-\frac{1}{q+1-p}}$$
by estimate $(\ref{est0})$, the result follows by integration. \qeda\medskip

An important consequence of the gradient estimate is the Harnack inequality.

\bprop{harnack} Assume $q>p-1$ and let $u\in C^1(\Gw)$ be a positive solution of $(\ref{A1})$ in $\Gw$. Then there exists a constant $C=C(n,p,q)>0$ such that for any $a\in \Gw$ and $R>0$ such that $\overline B_R(a)\subset \Gw$, there holds
\bel{harn1}
\max\{u(x):x\in B_{R/2}(a)\}\leq C\min\{u(x):x\in B_{R/2}(a)\}.
\ee
\es
\Proof We can assume $a=0$ in $\Gw$ and $R<d(0)=\dist (0,\prt\Gw)$. Then we write $(\ref{A1})$
$$-\Gd_p u+C(x)\abs{\nabla u}^{p-1}=0
$$
with $\abs{C(x)}=\abs{\nabla u}^{q+1-p}\leq c_{N,p,q}R^{-1}$ by $(\ref{est0})$. Set $u_R(y)=u(R y)$, then $u_R$ satisfies
$$-\Gd_p u_R+RC(Ry)\abs{\nabla u_R}^{p-1}=0\qquad\text{in }B_1.
$$
Since $RC(Ry)$ is bounded in $B_1$, we can apply Serrin's results (see \cite{Ser0}) and obtain
\bel{harn2}
\max\{u_R(y):y\in B_{1/2}(0)\}\leq C\min\{u(y):y\in B_{1/2}(0)\}.
\ee
Then $(\ref{harn1})$ follows.\qeda\medskip

The following Liouville result which improves a previous one due to Farina and Serrin \cite[Th 7]{FaSe}, is a direct consequence of the gradient estimate.


\bcor{Liouville} Assume $q>p-1>0$. Then any signed solution of $(\ref{A1})$ in $\BBR^N$ is a constant. \es
\Proof We apply of $(\ref{est0})$ in $B_R(a)$ for any $R>0$ and $a\in\BBR^N$ and let $R\to\infty$. \qeda



\mysection{Singularities in a domain}

\subsection{Radial solutions}
If $u$ is a radial function, we put $u(x)=u(\abs x)=u(r)$, with $r=\abs x$. If $u$ is a radial solution of $(\ref{A1})$ in $B_1^*:=B_1\setminus\{0\}$, it satisfies
\bel{rad1}
\left(\abs{u'}^{p-2}u'\right)'+\frac{N-1}{r}\abs{u'}^{p-2}u'-\abs{u'}^q=0
\ee
in $(0,1)$. We suppose $q<p$, then $p-1<q_c<p\leq N$. We set

\bel{rad2}b=\myfrac{N(p-1)-(N-1)q}{q+1-p}=\frac{(N-1)(q_c-q)}{q+1-p}.
\ee

 The next result provides the classification of radial solutions according to their sign near $0$.

\bprop{rad} Let $u$ be a nontrivial solution of $(\ref{rad1})$, then 
\bel{rad2'}
u'(r)=\left\{\BA {lll}
-r^{\frac{1-N}{p-1}}\left(b^{-1}r^{\frac{q+1-p}{p-1}b}+K\right)^{\frac{1}{p-1-q}}\quad&\text{if }\;q\neq q_c,\\
-r^{\frac{1-N}{p-1}}\left(\abs{\ln r^{N-1}}+K\right)^{\frac{1-N}{p-1}}\quad&\text{if }\;q= q_c.
\EA\right.\ee
As a consequence there holds\smallskip

\noindent {1- If $u$ is positive near $0$} and  $p-1<q<q_c$, \smallskip

\noindent (i) either there exists $k>0$ such that 
\bel{rad3}u(r)=\left\{\BA {ll}k\frac{N-p}{p-1}r^{\frac{p-N}{p-1}}+O(r^{\frac{q+1-N}{p-1}b}\vee 1)\quad&\text{if } p<N,\\[2mm]
-k\ln r+O(1)\quad&\text{if } p=N,
\EA\right.
\ee
and $u$ is a radial solution of
\bel{rad4}
-\Gd_pu+\abs{\nabla u}^q=c_{N,p}k^{p-1}\gd_0\qquad\text{in }\CD'(B_1),
\ee
(ii) or 
\bel{rad5}
u(r)=\gl_{N,p,q}r^{-\gb_q}+M,
\ee
where
\bel{rad6}
\gl_{N,p,q}=\gb_q^{-1}b^{\frac{1}{q+1-p}}.
\ee
\smallskip

\noindent   If $u$ is positive near $0$ and $q\geq q_c$, then $u$ is constant.\smallskip

\noindent {2-  If $u$ is negative near $0$}: then for  $p-1<q<q_c$,  there exists $k<0$ such that 
\bel{rad7}u(r)=\left\{\BA {ll}k\frac{N-p}{p-1}r^{\frac{p-N}{p-1}}+O(r^{\frac{q+1-N}{p-1}b}\vee 1)\quad&\text{if } p<N,\\[2mm]
-k\ln r+O(1)\quad&\text{if } p=N,
\EA\right.
\ee
and $u$ is a radial solution of
\bel{rad8-1}
\Gd_pu-\abs{\nabla u}^q=-c_{N,p}(-k)^{p-1}\gd_0\qquad\text{in }\CD'(B_1).
\ee
\smallskip

\noindent   If $q=q_c$, then\smallskip
\bel{rad9}u(r)=\left\{\BA {ll}-\gn_{N,p}r^{\frac{p-N}{p-1}}\left(-\ln r\right)^{-\frac{p-1}{N-1}}(1+o(1))\quad&\text{if } p<N\\[2mm]
-\gn_{N}\ln(-\ln r)(1+o(1))\quad&\text{if } p=N
\EA\right.
\ee
for some for some constant $\gn_{N,p},\,\gn_{N}>0$.\smallskip

\noindent   If $q>q_c$, 
\bel{rad10}
u(r)=-\gl_{N,p,q}r^{-\gb_q}+M.\ee
\es
\Proof
We set
\bel{rad11}
w(r)=r^{N-1}\abs{u'}^{p-2}u',
\ee
then
\bel{rad12}
w'(r)=r^{-\frac{(q+1-p)(N-1)}{p-1}}\abs{w}^{\frac{q}{p-1}}
\ee
Thus
\bel{rad13}
-\abs{w}^{-\frac{q}{p-1}}w=\left\{\BA{ll}
b^{-1}r^{\frac{q+1-p}{p-1}b}+K\qquad&\text{if }\,q\neq q_c\\[2mm]
\ln\left(Kr^{\frac{1}{N-1}}\right)\qquad&\text{if }\,q=q_c
\EA\right.
\ee
for some $K$. \smallskip	

\noindent{\it 1-Case $p-1<q<q_c$}, then $b>0$. If  $K>0$ then $w'$ and $u'$ are negative and 
\bel{rad13'}
u'(r)=-r^{\frac{1-N}{p-1}}\left[b^{-1}r^{\frac{q+1-p}{p-1}b}+K\right]^{-\frac{1}{q+1-p}}=-k'r^{\frac{p-N}{p-1}}+O(r^{\frac{q+2-p-N}{p-1}}).
\ee
Integrating again, we get $(\ref{rad3})$
From the asymptotic of $u'(r)$ we derive that  $u$ is a radial solution of $(\ref{rad4})$.
If $K=0$, then $u'(r)=-r^{-\frac{N-1+b}{p-1}}b^{\frac{1}{q+1-p}}$ and we get $(\ref{rad5})$, $(\ref{rad6})$.
This is the explicit particular solution.\smallskip

\noindent If $K=-\tilde K<0$, then $w'>0$ near $0$. We set $\tilde w=-w$, $\tilde u=-u$ and $\tilde w(r)=r^{N-1}\abs{\tilde u'}^{p-2}\tilde u'$. Thus, 
\bel{rad8-2}\tilde u(r)=\tilde k''r^{\frac{p-N}{p-1}}+O(r^{\frac{q+1-N}{p-1}b}\vee 1) \text{ or }u(r)=-\tilde k''\ln r+O(1),
\ee
according $p<N$ or $p=N$, and $\tilde u$ satisfies
\bel{rad8-3}
-\Gd_p\tilde u-\abs{\nabla u}^q=c_{N,p}\tilde k'\gd_0\qquad\text{in }\CD'(B_1).
\ee\smallskip	

\noindent {\it 2-Case $q\geq q_c$}. Then $b\leq 0$. If  $q>q_c$ (equivalently $b< 0$), $(\ref{rad5})$ implies 
\bel{rad9-1}
u'(r)=r^{\frac{1-N}{p-1}}\left[-b^{-1}r^{\frac{q+1-p}{p-1}b}+K\right]^{-\frac{1}{q+1-p}}=(-b)^{\frac{1}{q+1-p}}r^{-\frac{1}{q+1-p}}
(1+o(1)),\ee
then
\bel{rad10-1}
u(r)=-\gl_{N,p,q}r^{-\gb_q}
(1+o(1)).\ee
If $q=q_c$,
\bel{rad11-1}
u'(r)=r^{\frac{1-N}{p-1}}\left[(1-N)^{-1}\ln r+K\right]^{-\frac{p-1}{N-1}}
(1+o(1)),\ee
and, either $p<N$ and
\bel{rad12-1}
u(r)=-\gn_{N,p}r^{\frac{p-N}{p-1}}\left(-\ln r\right)^{-\frac{p-1}{N-1}}(1+o(1)),\ee
or $p=N$ and
\bel{rad12-2}
u(r)=-\gn_{N}\ln(-\ln r)(1+o(1)),\ee
for some constant $\gn_{N,p},\,\gn_{N}>0$.\qeda

\bprop{existrad}Assume $1<p\leq N$ and $p-1<q<q_c$, then for any $k>0$ there exists a unique positive solution $u=u_k$ of $(\ref{rad1})$ in 
$(0,1)$ vanishing for $r=1$ satisfying 
\bel{rad13-1}
\displaystyle\lim_{r\to 0}\frac{u_k(r)}{\gm_p(r)}=k.
\ee
When $k\to\infty$, $u_k\uparrow u_\infty$ which is a solution of $(\ref{rad1})$ in 
$(0,1)$ vanishing on $\prt B_1$ satisfying 
\bel{rad14}
\lim_{r\to 0}r^{\gb_q}u_\infty(r)=\gl_{N,p,q}.
\ee
\es
\Proof Using $(\ref{rad13'})$ we see that $K$ is completely determined by $K=k^{p-1-q}$ and $u_k$ by
\bel{rad15}
u_k(r)=\myint{r}{1}s^{\frac{1-N}{p-1}}\left[b^{-1}s^{\frac{q+1-p}{p-1}b}+k^{p-1-q}\right]^{-\frac{1}{q+1-p}}ds.
\ee
Conversely, asymptotic expansion in $(\ref{rad15})$ yields to $(\ref{rad13})$. The unique characterization of $K$ yields to uniqueness although uniqueness is also a consequence of the maximum principle as we will see it in the non radial case. Clearly the function $u_k$ defined by $(\ref{rad15})$ is increasing and $u_\infty=\lim_{k\to\infty}u_k$ satisfies
\bel{rad16}
u_\infty(r)=\myint{r}{1}s^{\frac{1-N}{p-1}}\left[b^{-1}s^{\frac{q+1-p}{p-1}b}\right]^{-\frac{1}{q+1-p}}ds=\gl_{N,p,q}(r^{-\gb_q}-1).
\ee
\bprop{globrad}Assume $1<p\leq N$ and $p-1<q<q_c$. If $u$ is a nonnegative radial solution of $(\ref{rad1})$ in 
$(0,\infty)$. Then\smallskip

\noindent (i) either $u(r)\equiv M$ for some  $M\geq 0$, or

\noindent (ii) there exist $k> 0$ and $M\geq 0$ such that
\bel{rad17}
u(r)=u_{k,M}(r):=\myint{r}{\infty}s^{\frac{1-N}{p-1}}\left[b^{-1}s^{\frac{q+1-p}{p-1}b}+k^{p-1-q}\right]^{-\frac{1}{q+1-p}}ds +M,
\ee

\noindent (ii) or there exists $M\geq 0$ such that
\bel{rad18}
u(r)=u_{\infty,M}(r):=\gl_{N,p,q}r^{-\gb_q}+M
\ee
\es
\Proof From identity $(\ref{rad13'})$, valid for any {\it nonconstant} solution $u$, we see that for a global positive solution we must have $K\geq 0$. If $K=0$ then $u=u_\infty$ defined by $(\ref{rad17})$. If $K>0$, then $u'\in L^1(1,\infty)$, thus
$u(\infty)=\lim_{s\to\infty} u(s)$ exists and
\bel{rad19}
u(r)=u(\infty)+\myint{r}{\infty}s^{\frac{1-N}{p-1}}\left[b^{-1}s^{\frac{q+1-p}{p-1}b}+K\right]^{-\frac{1}{q+1-p}}ds,
\ee
and $K=k^{p-1-q}$ in order to have $(\ref{rad17})$.\qeda

\subsection{Removable singularities}
\subsubsection{Removable singularities of renormalized solutions with right-hand side measures}

In this section $\Gw$ is any domain of $\BBR^N$. We denote by $\frak M(\Gw)$ the set of Radon measures in $\Gw$ and we study a more general equation than $(\ref{A1})$ 
\bel{A1+}
-\Gd_pu+|\nabla u|^q=\gm,
\ee
where $\gm\in \frak M(\Gw)$. For any $r>1$, the $C_{1,r}$ capacity is defined by
$$C_{1,r}(K)=\inf\{\norm{\psi}^r_{W^{1,r}}:\psi\in C_c^\infty(\BBR^N), \chi_K\leq\psi\leq 1\}
$$ 
for any compact subset $K$ of $\BBR^N$, and extended to capacitable sets by the classical method. We set
$$\frak M^r(\Gw)=\{\gm\in \frak M(\Gw):\gm (K)=0\;\forall\,K\subset\Gw,\; K\text{ compact}\text{ s.t. }C_{1,r}(K)=0\}.
$$ 
We recall that any measure $\gm$ in $\Gw$ can be decomposed in a unique way as
\bel{Decomp}
\gm=\gm_0+\gm_s^+-\gm_s^-
\ee
where $\gm_0\in \frak M^r(\Gw)$ and $\gm_s^\pm$ are nonnegative measures concentrated on sets with zero $C_{1,r}$ capacity.

In order to study equation $(\ref{A1+})$ it is natural to introduce other notions of  solutions than the strong ones. We use the notion of locally renormalized solutions introduced in \cite{B-V}, which gives a local version of the notion of renormalized solutions very much used in \cite{DMOP}, \cite {MuPor}, \cite{Mae}. 

For $k>0$ and $s\in\BBR$, we define the truncation $T_k(s)=\max\{-k,\min\{k,s\}\}$. If $u$ is measurable and finite a.e. and if $T_k(u)\in W^{1,p}_{loc}(\Omega)$, we define the gradient a.e. of $u$ by $\nabla T_k(u)=\chi_{\abs u\leq k}\nabla u$, for any $k>0$.
We denote by $q_*$ the conjugate exponent of $\frac{q}{p-1}$
\be\label{S1-1}
q_*= \frac{q}{q+1-p},
\ee
i.e. the conjugate of $q$ if $p=2$. 

\bdef{renorm} Let $u$ be a measurable and finite a.e. function in $\Gw$. Let $\gm=\gm_0+\gm_s^+-\gm_s^-\in \frak M(\Gw)$ with $\gm_0\in \frak M^p(\Gw)$ and $\gm_s^\pm$ singular and nonnegative as in $(\ref{Decomp})$.

\noindent 1- We say that $u$ is weak solution  of $(\ref{A1+})$ if $T_k(u)\in W^{1,p}_{loc}(\Omega)$ for any $k>0$, $\abs{\nabla u}^q\in L^1_{loc}(\Gw)$ and $(\ref{A1+})$ holds in the sense of distributions in $\Gw$.

\noindent 2- Assuming $\Gw$ is bounded, we say that $u$ is a renormalized (abridged R-solution) solution of $(\ref{A1})$ such that $u=0$ on $\prt\Gw$
if $T_k(u)\in W^{1,p}_{0}(\Omega)$ for any $k>0$, $\abs{\nabla u}^q\in L^{1}(\Omega)$ and
\be\label{S1-2}
\abs u^{p-1}\in L^\gs(\Omega),\;\forall\gs\in [1,\frac N{N-p})\text{ and }
\abs {\nabla u}^{p-1}\in L^\gt(\Omega),\;\forall\gt\in [1,\frac N{N-1}),
\ee
 and for any $h\in W^{1,\infty}(\BBR)$ such that $h'$ has compact support and any $\gf\in W^{1,m}(\Omega)$ for some $m>N$   such that 
 $h(u)\gf\in W_0^{1,m}(\Omega)$ there holds
 \be\label{S1-2+}\BA {lll}
\myint{\Gw}{}\left(|\nabla u|^{p-1}\nabla u.\nabla(h(u)\gf) +|\nabla u|^{q}h(u)\gf\right)dx\\[4mm]
\phantom{----------}
=\myint{\Gw}{}h(u)\gf d\gm_0+h(\infty)\myint{\Gw}{}\gf d\gm^+_s-h(-\infty)\myint{\Gw}{}\gf d\gm^-_s.
\EA\ee

\noindent 3- We say that $u$ is a local renormalized (abridged LR-solution)  of solution of $(\ref{A1})$ if $T_k(u)\in W^{1,p}_{loc}(\Omega)$ for any $k>0$, $\abs{\nabla u}^q\in L^{1}_{loc}(\Omega)$ and
\be\label{S1-2++}
\abs u^{p-1}\in L^\gs_{loc}(\Omega),\;\forall\gs\in [1,\frac N{N-p})\text{ and }
\abs {\nabla u}^{p-1}\in L^\gt_{loc}(\Omega),\;\forall\gt\in [1,\frac N{N-1})
\ee
 and for any $h\in W^{1,\infty}(\BBR)$ such that $h'$ has compact support and any
 $\gf\in W^{1,m}(\Omega)$ for some $m>N$ with compact support and  such that 
 $h(u)\gf\in W_0^{1,m}(\Omega)$ identity $(\ref{S1-2+})$ holds.

\es

\noindent\Remark If $q\geq 1$ and $u$ is a weak solution, then it satisfies $(\ref{S1-2+})$, see 
\cite[Lemmas 2.2, 2.3]{DMOP}.

\medskip

Our main removability result is the following.

\bth{remov} Assume $0<p-1<q\leq p$. If $F\subset\Gw$ is a relatively closed set such that $C_{1,q_*}(F)=0$, and $\gm\in \mathfrak M^{q^*}(\Gw)$.\smallskip

\noindent (i) Let $p-1<q\leq p$ and $u$ be a LR-solution of  $(\ref{A1+})$ in $\Gw\setminus F$. Then $u$ is a LR-solution of  $(\ref{A1+})$ in $\Gw$.
\smallskip

\noindent (ii) Let $q>p$ and $u$ be a weak solution of $(\ref{A1+})$ in $\Gw\setminus F$. Then $u$ is a  weak solution of  $(\ref{A1+})$ in $\Gw$.
\es
\Proof Notice that a set $F$ with $C_{1,q_*}(F)=0$, has zero measure; since $u$ is defined up to a set of zero measure, any extension of $u$ to $F$ is valid. Notice also that if $p-1<q<q_c$ then $W^{1,q_*}(\BBR^N)$ is imbedded into $C(\BBR^N)$, therefore only the empty set has zero $C_{1,q_*}$ capacity.\smallskip

\noindent (i) From our assumption $T_k(u)\in W^{1,p}_{loc}(\Gw\setminus\{F\})$ for any $k>0$,  $|u|^{p-1}\in L^\gs_{loc}(\Omega)$ for all $\gs\in [1,\frac{N}{N-p})$,
  $|\nabla u|^{p-1}\in L^\gt_{loc}(\Omega)$ for all $\gt\in [1,\frac{N}{N-1})$, and  $|\nabla u|^{q}\in L^1_{loc}(\Omega)$. Since $p\leq q_*$, for any compact $K\subset\Omega$, $C_{1,p}(F\cap K)=0$. Thus $T_k(u)\in W^{1,p}_{loc}(\Gw)$ by \cite[Th 2.44]{HKM}. Because $u$ is measurable and finite a.e. on $\Gw$, we can define $\nabla u$ a.e. in $\Omega$ by the formula $\nabla u=\nabla T_k(u)$ a.e. on the set $\{x\in\Omega:|u(x)|\leq k\}$.
  
  Let $\gz\in C_c^\infty(\Gw)$ with support in $\gw\subset\bar\gw\subset \Gw$, $\gz\geq 0$. Set $K_\gz=F\cap supp\,\gz$. Then $K_\gz$ is compact and $C_{1, q_*}(K_\gz)=0$, thus there exists $\gz_n\in C_0^\infty(\Gw)$ such that $0\leq\gz_n\leq 1$, $\gz_n=1$ in a neighborhood of $K_\gz	$ that we can assumed to be contained in $\gw$, such that $\gz_n\to 0$
  in $W^{1, q_*}(\BBR^N)$. It can also be assumed that $\gz_n(x)\to 0$ for all $x\in\BBR^N\setminus E$ where 
  $E$ is a Borel set such that  $C_{1, q_*}(E)=0$ (see e.g. \cite[Lemmas 2.1, 2.2]{BP}). Let $\xi_n=\gz(1-\gz_n)$. Since $u$ is a weak solution of $(\ref{A1+})$ in $\Gw\setminus F$, we can take $\xi_n^{q_*}$ as a test function and get
  $$\myint{\Gw}{}\left(|\nabla u|^q\xi_n^ {q_*}+ {q_*} \xi_n^{ {q_*}-1}|\nabla u|^{p-2}\nabla u.\nabla \xi_n\right) dx=\myint{\Gw}{}\xi_n^{q_*} d\gm.
  $$
By H\"older's inequality, for any $\eta>0$,
$$\BA {lll}\myint{\Gw}{}|\nabla u|^q\xi_n^ {q_{*}} dx\leq  {q_*} \myint{\Gw}{}\xi_n^{ {q_*}-1}|\nabla u|^{p-1}|\nabla \xi_n| dx\\[4mm]
\phantom{\myint{\Gw}{}|\nabla u|^q\xi_n^ {q_{*}} dx}
\leq ( {q_*} -1)\eta^{\frac{q}{p-1}}\myint{\Gw}{}|\nabla u|^q\xi_n^ {q_*} dx+\eta^{- {q_*}}\myint{\Gw}{}|\nabla \xi_n|^ {q_*} dx.
\EA$$
Hence, taking $\eta$ small enough,
$$\BA {lll}\myint{\Gw}{}|\nabla u|^q\xi_n^ {q_*} dx\leq c\left(\myint{\Gw}{}|\nabla \xi_n|^ {q_*} dx+\abs{\myint{\Gw}{}\xi_n^ {q_*}d\gm}\right)
\\[4mm]\phantom{\myint{\Gw}{}|\nabla u|^q\xi_n^ {q_*} dx}
\leq c\left(\myint{\Gw}{} \left(|\nabla \gz|^ {q_*} +|\nabla \gz_n|^ {q_*}\right) dx+
\myint{{\rm supp}(\gz)}{}d\abs\gm\right).
\EA$$
with $c=c(p,q)>0$. From Fatou's lemma, we get $|\nabla u|^q\gz^ {q_*}\in L^1(\Gw)$ and
\be\label{S1-7}
\myint{\Gw}{}|\nabla u|^q\gz^ {q_*} dx\leq c_\gz:=c\left(\myint{\Gw}{}|\nabla \gz|^ {q_*} dx+
\myint{{\rm supp}(\gz)}{}d\abs\gm\right).
\ee
Taking now $T_k(u)\xi_n^ {q_*}$ as test function, we obtain
$$\BA {lll}\myint{\Gw}{}|\nabla (T_k(u))|^p\xi_n^ {q_*} dx+
\myint{\Gw}{}|\nabla u|^qT_k(u)\xi_n^ {q_*} dx=-\myint{\Gw}{}T_k(u)|\nabla u|^{p-2}\nabla u.\nabla \xi_n^ {q_*} dx\\[4mm]
\phantom{----\myint{\Gw}{}|\nabla u|^qT_k(u)\xi_n^ {q_*} dx}
+\myint{\Gw}{}T_k(u)\xi_n^ {q_*} d\gm_0+k\left(\myint{\Gw}{}\xi_n^ {q_*} (d\gm_s^+-d\gm_s^-\right).
\EA$$
Then we deduce, from H\"older's inequality,
$$\BA {llll}
\myfrac 1k\abs{\myint{\Gw}{}T_k(u)|\nabla u|^{p-2}\nabla u.\nabla \xi_n^ {q_*} dx}
\leq  {q_*} \myint{\Gw}{}\left(\gz^{ {q_*}-1}|\nabla u|^{p-1}|\nabla\gz|+ 
\gz^{ {q_*}}|\nabla u|^{p-1}|\nabla\gz_n|\right) dx
\\[4mm]\phantom{\myfrac 1k\abs{\myint{\Gw}{}T_k(u)|\nabla u|^{p-2}\nabla u.\nabla \xi_n^ {q_*} dx}}
\leq (2 {q_*} -1)\myint{\Gw}{}\left(|\nabla u|^q\gz^ {q_*}+|\nabla \gz|^ {q_*}+ {q_*} \gz^ {q_*} |\nabla \gz_n|^ {q_*}\right) dx
\\[4mm]\phantom{\myfrac 1k\abs{\myint{\Gw}{}T_k(u)|\nabla u|^{p-2}\nabla u.\nabla \xi_n^ {q_*} dx}}
\leq 2 {q_*} c_\gz+\myint{\Gw}{}|\nabla \gz|^ {q_*} dx+o(1).
\EA
$$
Therefore, up to changing $c_\gz$ into another constant $c_\gz$ depending on $\gz$,
$$\myint{\Gw}{}|\nabla (T_k(u))|^p\xi_n^ {q_*} dx
\leq (k+1)c_\gz+o(1),
$$
and by Fatou's lemma,
\be\label{S1-8}
\myint{\Gw}{}|\nabla (T_k(u))|^p\gz ^ {q_*} dx
\leq (k+1)c_\gz.
\ee
By a variant of the results in \cite{BoGa},\cite{BoGa1} due to \cite {PPP} it follows that the regularity statements $(\ref{S1-2++})$ of \rdef{renorm} hold. \smallskip

Finally, we show that $u$ is a LR-solution in $\Gw$. Let $h\in W^{1,\infty}(\BBR)$ such $h'$ has compact support and let $\gf\in W^{1,m}(\Gw)$ with $m>N$ with compact support in $\Gw$, such that $h(u)\gf\in W^{1,p}(\Gw)$. Consider again $\gz$, $\gz_n$ and $\xi_n$ as above. Then $(1-\gz_n)\gf\in W^{1,m}(\Gw\setminus F)$ and $h(u)(1-\gz_n)\gf\in W^{1,p}(\Gw\setminus F)$ and has compact support in $\Gw\setminus F$. We can write
$$I_{1,n}+I_{2,n}+I_{3,n}+I_{4,n}=I_{5,n}+I^+_{6,n}+I^-_{6,n},
$$
where
$$
I_{1,n}=\myint{\Gw}{}|\nabla u|^{p}h'(u)(1-\gz_n)\gf  dx\,,\; I_{2,n}=-\myint{\Gw}{} h(u)\gf|\nabla u|^{p-2}\nabla u.\nabla \gz_n dx,
$$
$$I_{3,n}=\myint{\Gw}{} h(u)(1-\gz_n)|\nabla u|^{p-2}\nabla u.\nabla \gf dx
\,,\; I_{4,n}=\myint{\Gw}{}|\nabla u|^{q}(1-\gz_n)\gf  dx.
$$
$$I_{5,n}=\myint{\Gw}{}h(u)\gf(1-\gz_n)d\gm_0\,\text{ and }\; I^\pm_{6,n}=h(\pm\infty)\myint{\Gw}{}\gf(1-\gz_n)d\gm^\pm_s.
$$
We get $\lim_{n\to\infty} I_{1,n}=\myint{\Gw}{}|\nabla u|^{p}h'(u)\gf  dx$ since there exists some $a>0$, independent of $n$, such that 
$$I_{1,n}=\myint{\Gw}{}\abs{\nabla T_a(u)}^{p}h'( T_a(u))(1-\gz_n)\gf  dx.
$$
Furthermore $\lim_{n\to\infty}I_{2,n}=0$ since
$$\abs{\myint{\Gw}{} h(u)\gf|\nabla u|^{p-2}\nabla u.\nabla \gz_n dx}\leq \norm h_{L^\infty}\left(\myint{\Gw}{}\abs{\nabla u}^q\right)^{\frac{p-1}{q}}\norm{\nabla\gz_n}_{L^{q_*}}.
$$
Next $\lim_{n\to\infty} I_{3,n}=\myint{\Gw}{}h(u)|\nabla u|^{p-2}\nabla u.\nabla\gf  dx$ because $\nabla\gf\in L^m(\Gw)$ and $|\nabla u|^{p-1}\in L^\gt_{loc}(\Gw)$ for all $\gt\in [1,\frac{N}{N-1})$. From $(\ref{S1-7})$
$\lim_{n\to\infty}I_{4,n}=\myint{\Gw}{}h(u)\gf\abs{\nabla u}^qdx$. But $h(u)\gf\in L^1(\Gw,d\abs{\gm_0})$ by
\cite[Remark 2.26]{DMOP}. Since $\gz_n\to 0$ everywhere in $\BBR^N\setminus E$ and $\gm(E)=0$, it follows $\lim_{n\to\infty} I_{5,n}=\myint{\Gw}{}h(u)\gf  d\gm_0$. Clearly 
$\lim_{n\to\infty}I^\pm_{6,n}=h(\pm\infty)\myint{\Gw}{}\gf d\gm^\pm_s$. Hence $u$ is a LR solution in whole $\Gw$.\smallskip

\noindent (ii) Let $u$ be a weak solution in $\Gw\setminus F$. SInce $q>p$, $1<q_*<p$, hence $u\in W^{1,q_*}_{loc}(\Gw\setminus F)=W^{1,q_*}_{loc}(\Gw)$ and $\abs{\nabla u}\in L^1_{loc}(\Gw)$. Let 
$\gz\in C^\infty_c(\Gw)$ and $\gz_n$ and $\xi_n$ as in (i). We obtain again $|\nabla u|^{q}\gz^{q_*}\in L^1(\Gw)$. Hence $\nabla u\in L^q_{loc}(\Gw)$. Next we take $\xi_n$ as a test function in equation $(\ref{A1+})$ in $\Gw\setminus F$. We obtain $J_{1,n}+J_{2,n}+J_{3,n}=J_{4,n}$ with
$$J_{1,n}=\myint{\Gw}{}(1-\gz_n)\abs{\nabla u}^{p-2}\nabla u.\nabla \gz dx\,,\; 
J_{2,n}=-\myint{\Gw}{}\gz\abs{\nabla u}^{p-2}\nabla u.\nabla \gz_n dx,
$$
$$J_{3,n}=\myint{\Gw}{}\abs{\nabla u}^{q}\gz(1-\gz_n) dx\,,\;J_{4,n}=\myint{\Gw}{}\gz(1-\gz_n) d\gm.
$$
We can let $n\to\infty$ in $J_{1,n}$ and $J_{3,n}$ using the dominated convergence theorem and the fact that $\nabla u\in L^q_{loc}(\Gw)$ and $q>p-1$. Furthermore $\lim_{n\to\infty}J_{2,n}=0$ because 
$|\nabla u|^{p-1}\in L^{\frac{q}{p-1}}_{loc}(\Gw)$ and $|\nabla\gz_n|\to 0$ in $L^{q_*}(\Gw)$. Since 
$J_{4,n}\to \myint{\Gw}{}\gz d\gm$ as above, it follows that $u$ is a weak solution in $\Gw$.
\qeda

\subsubsection{Regularity results}
The natural question concerning LR-solutions obtained in \rth{remov} is their regularity. It is noticeable that the results are very different according to whether we consider nonnegative or signed solutions.   Here we give some regularity properties of solutions of $(\ref{A1})$. We first consider nonnegative solutions of $(\ref{A1})$.

\bth{Reg3} Let $p-1<q$, $N\geq 2$ and $u$ is a nonnegative LR-solution of $(\ref{A1})$. Then $u\in L^\infty_{loc}(\Gw)\cap W^{1,p}_{loc}(\Gw)$. As a consequence, if $q\leq p$, $u\in C^{1,\ga}(\Gw)$ for some $\ga\in (0,1)$.
\es
\Proof Since $-\Gd_p\leq 0$, then $u\in L^\infty_{loc}(\Gw)$ by a recent argument due to Kilpelainen and Kuusi \cite{KKTK} and $u$ satisfies the weak Harnack inequality
$$\sup_{B_\gr(x_0)}\leq \gr^{-N}\left(\myint{B_{2\gr}(x_0)}{}u^ {q_*} dx\right)^{\frac{1}{ {q_*}}}
$$
with $C=C(N,p, {q_*})$. Then $u$ coincides with $T_k(u)$ in any ball $B_\gr(x_0)$ such that $\overline B_{2\gr}(x_0)\subset\Gw$, for $k$ large enough. Thus 
$u\in W^{1,p}_{loc}(\Gw)$. If $q\leq p$, it follows by Tolksdorff's result \cite{To} that $u\in C^{1,\ga}(\Gw)$.\qeda\medskip

When we  deal with signed solutions of $(\ref{A1})$, there is another critical value
involved when $q\leq p$,
\bel{R0}
\tilde q=p-1+\frac{p}{N}.
\ee
Observe that $q_c<\tilde q<p$ if $1<p<N$ and $q_c=\tilde q=N$ if $p=N$. For simplicity we consider solutions of
\be\label{R1X}\BA {lll}
-\Gd_pu+|\nabla u|^q=0\qquad&\text{in }\Gw\\
\phantom{-\Gd_p+|\nabla u|^q}
u=0\qquad&\text{on }\prt\Gw,
\EA\ee
and we first recall some local estimates of the gradient of renormalized solutions. 
\blemma{Reg1} Assume $\Gw$ is a bounded $C^2$ domain. Let $u$ be a renormalized solution of the problem
\be\label{R1}\BA {lll}
-\Gd_pu=f\qquad&\text{in }\Gw\\
\phantom{-\Gd_p}
u=0\qquad&\text{on }\prt\Gw
\EA\ee
where $f\in L^m(\Gw)$ with $1<m<N$ and set $\bar m=\frac{Np}{Np-N-p}=\frac{p}{\tilde q}$, where $\tilde q$ is defined in $(\ref{S1-1})$.\smallskip

\noindent (i) If $m>\frac{N}{p}$, then $u\in L^\infty(\Gw)$. If $m=\frac{N}{p}$,	then 
$u\in L^k(\Gw)$ for $1\leq k<\infty$. If  $m<\frac{N}{p}$, then $|u|^{p-1}\in L^k(\Gw)$ with $k=\frac{Nm}{N-mp}$.\smallskip

\noindent (ii) $\nabla u^{p-1}\in L^{m^*}(\Gw)$ with $m^*=\frac{Nm}{N-m}$. Furthermore, if $\bar m\leq m$, then $u\in W^{1,p}_0(\Gw)$.
\es
\Proof The estimates in the case $m<\bar m$ are obtained in \cite{B-VH} following \cite{BoGa1} and \cite{KL}, by using for test functions $\gf_{\gb,\ge}(T_k(u))$ where
$$\gf_{\gb,\ge}(w)=\myint{0}{w}(\ge+|t|)^{-\gb}dt
$$
for some $\gb<1$. In the case $m\geq \bar m$ and $1<p<N$ there holds $L^m(\Gw)\subset W^{-1,p'}(\Gw)$, thus $u$ is a variational solution in $W^{1,p}_0(\Gw)$. In the case $m=\bar m$, then $m^*=p'$ and the  conclusion follows. Next, if $m> \bar m$ or equivalently $m^*>p'$, then for any $\gs>p$ and 
$F\in (L^\gs(\Omega))^N$, there exists a unique $w\in W^{1,\gs}_0(\Gw)$, weak solution of 
 \be\label{R1-1}\BA {lll}
-\Gd_pw=div (|F|^{p-2}F)\qquad&\text{in }\Gw\\
\phantom{-\Gd_p}
w=0\qquad&\text{on }\prt\Gw,
\EA\ee
see \cite{Iv}, \cite{KZ}, \cite{KZ1}. Let $v$ be the unique solution in $W^{1,1}_0(\Gw)$
\be\label{R1-2}\BA {lll}
-\Gd v=f\qquad&\text{in }\Gw\\
\phantom{-\Gd}
v=0\qquad&\text{on }\prt\Gw.
\EA\ee
From the classical $L^p$-theory, $v\in W^{2,m}(\Gw)$, thus $\nabla v\in L^{m^*}(\Gw)$. Let $F$ be defined by $|F|^{p-2}F=\nabla v$. Then $F\in (L^\gs(\Omega))^N $ with $\gs=(p-1)m^*>p$. Then 
$$-\Gd_pw=-\Gd v=f.
$$
Thus $w=u$. This implies $u\in W^{1,\gs}_0(\Gw)$ and therefore $|\nabla u|^{p-1}\in L^{m^*}(\Gw)$.\qeda\medskip

Our first result is valid without any sign assumption on the solution.

\bth{Reg2} Assume $\Gw$ is a bounded $C^2$ domain. Let $p-1<q<\tilde q$, $N\geq 2$ and $u$ be a renormalized solution of problem $(\ref{R1X})$, such that 
\be\label{R1-3}\BA {lll}
|\nabla u|^q\in L^{m_0}(\Gw)\quad\text{for some }m_0>\max\{1,\frac{N(q+1-p)}q\}.
\EA\ee
Then $u\in C^{1,\ga}(\bar\Gw)$ for some $\ga\in (0,1)$. In particular $(\ref{R1-3})$ is satisfied if $q<q_c$, or if $q_c\leq q<\tilde q$ and $u\in W^{1,p}_0(\Gw)$.
\es
\Proof Set $f=-|\nabla u|^q\in L^{m_0}(\Omega)$. If $m_0\geq N$, $f\in L^{N-\gd}(\Gw)$ for any $\gd\in (0,N-1]$. Then $f\in L^{m_1}(\Omega)$ with $m_1=\frac{(p-1)m^*_0}{q}$. Note that $\frac{m_1}{m_0}=\frac{N(p-1)}{qN-p}>1$ since $q<\tilde q$. By induction, starting from $m_1$, we can defined $m_n$ as long as it is smaller than $N$ by $m_{n}=\frac{(p-1)m^*_{n-1}}{q}$, and we find $m_n<m_{n+1}$. If $m_n<N$ for any $n\in\BBN$, the sequence $\{m_n\}$
would converge to $L=\frac{N(q+1-p)}{q}$, which is impossible since we have assumed $m_0>L$. Therefore there exists some $n_0$ such that $m_{n_0}\geq N$. If $m_{n_0}=N$, (or if $m_{0}=N$ we can modify it so that $m_{n_0}<N$, but $m_{n_0+1}>N$. Then we conclude as above. 

If $q<q_c$, then $|\nabla u|^{p-1}\in L^{\frac{N(1-\gd)}{N-1}}(\Gw)$ for $\gd>0$ small enough. Then we can choose $m_0$ such that 
$\max\{1,\frac{N(q+1-p)}q\}<m_0<\frac{N}{N-1}.$ If $q_c\leq q<\tilde q$ and 
$u\in W^{1,p}_0(\Gw)$, we choose $m_0=\frac{p}{q}$.\qeda\medskip

\noindent\Remark The result which holds without sign assumption on $u$ is sharp. Indeed, the function $v=V+\tilde\gl_{N,p,q}$ defined above does not satisfy assumption $(\ref{R1-3})$, since $|\nabla u|^q\in L^{m}(\Gw)$  if and only if 
$m<\frac{N}{q(\gb_q+1)}=\frac{N(q+1-p)}{q}$. \medskip


\subsection{Classification of isolated singularities}
\subsubsection{Positive solutions}
The next result provides the complete classifications of isolated singularities of nonnegative solutions of $(\ref{A1})$. We suppose that $\Gw$ is an open subset of $\BBR^N$ containing $0$ and set $\Gw^*=\Gw\setminus\{0\}$. Without loss of generality, we can suppose that $\Gw\supset \overline B_1$ and we also recall that $B_1^*=B_1\setminus\{0\}$. We recall that the fundamental solution of the $p$-Laplacian is defined in $\BBR^N_*$ by
 \bel{abs-0}
 \gm_p(x)=\left\{\BA {ll}\abs x^{\frac{p-N}{p-1}}\quad&\text{if }1<p<N\\
 -\ln\abs x\quad&\text{if }p=N,
 \EA\right.
 \ee
and it satisfies
 \bel{abs-00}
-\Gd_p \gm_p=c_{N,p}\gd_0\quad\text{in }\CD'(\Omega).
 \ee
\bth{class} Let $p-1<q<q_c$ and $1<p\leq N$. If $u\in C^1(\Gw^*)$ is a nonnegative solution of $(\ref{A1})$ in $\Gw^*$, then we have the following alternative.\smallskip

\noindent (i) Either there exists $k\geq 0$ such that

\be\label{cl1}
\lim_{x\to 0}\myfrac{u(x)}{\gm_p(x)}=k
\ee
and $u$ satisfies
\be\label{cl2}
-\Gd_pv+\abs{\nabla v}^q=c_{N,p}k^{p-1}\gd_0\quad \text{in }\CD'(\Gw).
\ee
\smallskip

\noindent (ii) Or
\be\label{cl3}
\lim_{x\to 0}\abs x^{\gb_q}u(x)=\gl_{N,p,q},
\ee
where $\gb_q$ and  $\gl_{N,p,q}$ are defined in $(\ref{I-7})$. \smallskip

Furthermore, if $\Gw$ is bounded,  a nonnegative solution $u$ in $C(\overline\Gw\setminus\{0\})$ is uniquely determined by its data on $\prt\Gw$ and its behaviour $(\ref{cl1})$ or $(\ref{cl3})$ at $0$.
\es

We need several lemmas for proving this theorem. The method developed below for obtaining point wise estimates of the derivatives is an adaptation of a technique introduced in \cite{FV}.
\blemma{classlem0}Assume $p$, $q$ are as in \rth{class} and $\gf:(0,1]\mapsto\BBR_+$ is a continuous decreasing function such that $\gf(2r)\leq a\gf(r)$ and $r^{\frac{p-q}{q+1-p}}\gf(r)\leq c$ for some $a,c>0$ and any $r>0$.
If $u$ is a solution of $(\ref{A1})$ in $B_1^*$ such that 
\be\label{cl01}
\abs{u(x)}\leq \gf(\abs x)\qquad\forall x\in B_1^*.
\ee
Then there exists $C>0$ and $\ga\in (0,1)$, both depending on $N$, $p$, $q$, such that
\be\label{cl02}
\abs{\nabla u(x)}\leq C\gf(\abs x)\abs x^{-1}\qquad\forall x\in B_\frac{1}{2}\setminus\{0\}.
\ee
\be\label{cl03}
\abs{\nabla u(x)-\nabla u(x')}\leq C\gf(\abs x)\abs x^{-1-\ga}\abs{x-x'}^\ga\qquad\forall x,x'\text{ s.t. }0<\abs x\leq \abs {x'}\leq \frac{1}{2}.
\ee
\es
\Proof Define $\Gamma:=\{y\in\BBR^N:1<\abs y<7\}$ and $\Gamma'=\{y\in\BBR^N:2\leq \abs y\leq 6\}$. For $0<\abs x<\frac{1}{2}$ there exists $\ell\in (0,\frac{1}{4})$ such that $2\ell\leq\abs x\leq 3\ell$. We set 
$$u_\ell(y)=\frac{1}{\gf(\ell)}u(\ell y).$$
Then the equation
$$-\Gd_pu_\ell+\ell^{p-q}(\gf(\ell))^{q+1-p}\abs{\nabla u_\ell}^q=0
$$
holds in $\Gamma$. Because of $(\ref{cl01})$ and the fact that $\gf$ is decreasing, $u_\ell(y)\leq 1$ on $\Gamma$. Since $\ell^{p-q}(\gf(\ell))^{q+1-p}$ remains bounded for $\ell\in (0,1]$, we can apply Tolksdorff's a priori estimates \cite{To} and derive that
\be\label{cl04}
\abs{\nabla u_\ell(y)}\leq C\qquad\forall y\in \Gamma^*,
\ee
\be\label{cl05}
\abs{\nabla u_\ell(y)-\nabla u_\ell(y')}\leq C\abs{y-y'}^\ga \qquad\forall y,y'\in \Gamma^*,
\ee
for some $C=C\left(N,p,q,\norm {u_\ell}_{L^\infty(\Gamma)}\right)$ and $\ga\in (0,1)$. Putting $x=\ell y$, $x'=\ell y'$ where $x$, $x'$ are such that $0<\abs x\leq \abs {x'}\leq 2|x|\leq 1$ we have $|y'|=\frac{|x'|}{\beta}\leq \frac{2|x|}{\beta}\leq 6$ and thus
$$\abs{\nabla u(x)-\nabla u(x')}\leq C\gf(\ell)\ell^{-1-\ga}\abs{x-x'}^\ga\leq C\gf(\abs x)\abs x^{-1-\ga}\abs{x-x'}^\ga.
$$
If  $|x'|>2|x|$ we have
$$\abs{\nabla u(x)-\nabla u(x')}\leq C\left(\frac{\gf(|x|)}{|x|}+\frac{\gf(|x'|)}{|x|}\right)\leq \frac{2C\gf(|x|)}{|x|}
\leq \frac{2C\gf(|x|)}{|x|^{1+\ga}}|x-x'|^\ga.
$$
\qeda

\blemma{classlem1} Assume $p$, $q$ are as in \rth{class}. Let $u$ be a nonnegative solution of $(\ref{A1})$ in $\Gw^*$ such that 
\be\label{cl4}
\liminf_{x\to 0}\myfrac{u(x)}{\gm_p(x)}<\infty.
\ee
Then there exists $k\geq 0$ such that $(\ref{cl1})$ and $(\ref{cl2})$ hold.
\es
\Proof Let $y\in B_{\frac{1}{2}}$. By \rprop{harnack} there exists $C=C(N,p,q)>0$ such that
$$\max_{\abs{z-y}\leq \frac{\abs y}{4}}u(z)\leq C\min_{\abs{z-y}\leq \frac{\abs y}{4}}u(z)
$$
By a simple 2-D geometric construction and the help of trigonometric estimates it is easy to check that  if $\abs{y'}=\abs{y}$ there exist a chain of at most 7 balls $B_{\frac{\abs y}{4}}(y_j)$ with center  $y_j$ on $\{z:\abs z=\abs y\}$ such that 
$y_1=y$, $y_{13}=y'$ and $B_{\frac{\abs y}{4}}(y_j)\cap B_{\frac{\abs y}{4}}(y_{j+1})\neq\emptyset$. This implies
$u(y)\leq C^{7}u(y')$ and, since $\gm_p$ is radial (and we note $\gm_p(x)=\gm_p(|x|)$), 
\be\label{cl5}
\limsup_{x\to 0}\myfrac{u(x)}{\gm_p(x)}=k<\infty.
\ee
If $k=0$, then $u\leq \ge\gm_p+M$ for any $\ge>0$, where $M=\max\{u(x):\abs x=1\}$, by the comparison principle. Thus 
$u$ remains bounded near $0$ and therefore the singularity is removable. Next we assume $k>0$, thus $u\leq k(\gm_p-\gm_p(1))+M$ by applying the comparison principle in $\overline B_1\setminus\{0\}$.\smallskip

\noindent Up to changing $B_1$ into $B_{r_0}a$ for some $r_0\in (0,1)$ it implies $ u(x)\leq m\gm_p(x)$ for some $m\geq k$.
 Since $q\leq q_c$, $(\gm_p(r))^{q+1-p}r^{p-q}\leq c$, it follows from \rlemma{classlem0} that
\be\label{cl7}
\abs{\nabla u(x)}\leq C\gm_p(\abs x)\abs x^{-1}\qquad \forall x\in B_{\frac{r_0}{2}}\setminus\{0\}
\ee
and
\be\label{cl8}
\abs{\nabla u(x)-\nabla u(x')}\leq C\gm_p(\abs x)\abs x^{-1-\ga}\abs{x-x'}^{\ga}\quad\forall x,x'\text{ s.t. }0<
\abs x\leq\abs {x'}\leq\frac{r_0}{2}.\ee
If we define 
\be\label{cl9}u_r(y)=\frac{u(ry)}{\gm_p(r)}\qquad\forall y\in B_{\frac{r_0}{r}},
\ee
 it satisfies
\be\label{cl10}
-\Gd_pu_r+(\gm_p(r))^{q+1-p}r^{p-q}\abs{\nabla u_r}^q=0
\ee
in $B_{\frac{r_0}{r}}$ and the following estimates:
\be\label{cl10'}
0\leq u_r(y)\leq m\frac{\gm_p(r\abs y)}{\gm_p(r)}\qquad \forall y\in B_{\frac{r_0}{r}}\setminus\{0\},
\ee
\be\label{cl11}
\abs{\nabla u_r(y)}\leq C\frac{\gm_p(r\abs y)}{\gm_p(r)}\abs y^{-1}\qquad \forall y\in B_{\frac{r_0}{2r}}\setminus\{0\},
\ee
and
\be\label{cl12}
\abs{\nabla u_r(y)-\nabla u_r(y')}\leq C\frac{\gm_p(r\abs y)}{\gm_p(r)}\abs y^{-1-\ga}\abs{y-y'}^{\ga}\quad\forall y,y'\text{ s.t. }0<
\abs y\leq\abs {y'}\leq\frac{1}{2r_0}.\ee
Let $0<a<b$. If we assume that $0<a\leq\abs y\leq b$, then $\frac{\gm_p(r\abs y)}{\gm_p(r)}$ remains bounded independently of $r\in (0,1]$ and the set of functions $\{u_r\}_{0<r<1}$ is relatively sequentially compact in the $C^{1}$ topology of $\overline B_{b}\setminus B_{a}$. There exist a sequences $\{r_n\}$ converging to $0$ and a function $v\in C^{1}(\overline B_{b}\setminus B_{a})$ such that $u_{r_n}\to v$ in $C^{1}(\overline B_{b}\setminus B_{a}))$. Since 
$(\gm_p(r_n))^{q+1-p}r_n^{p-q}\to 0$ as $q<q_c$,  $v$ is p-harmonic in $B_{b}\setminus \overline B_{a}$ and nonnegative. Notice that 
$a$ and $b$ are arbitrary, therefore, using Cantor diagonal process, we can assume that $v$ is defined in $\BBR^N_*$ and $u_{r_n}\to v$ in the $C^1$-loc topology of $\BBR^N_*$. If $p=N$, the positivity of $v$ implies  that $v$ is a constant \cite[Corollary 2.2]{KV}, say $\gth$. If $1<p<N$, there holds, by  \cite[Theorem 2.2]{KV} and 
$(\ref{cl10'})$ 
\be\label{cl13}
v(y)=\gth\gm_p(y)+\gs\leq m\gm_p(y)\qquad\forall y\in\BBR^N_*,
\ee
for some $\gth,\gs\geq 0$, thus $\gs=0$. In order to make $\gth$ precise we set
$$\gamma(r)=\sup_{\abs x=r}\frac{u(x)}{\gm_p(x)},
$$
then $u(x)\leq \gamma (r)\gm_p(x)$ in $B_{r_0}\setminus B_r$. This implies in particular that, for $r<s<1$, 
$u(x)\leq \gamma (r)\gm_p(x)$ for any $x$ such that $\abs x=s$ and finally
\be\label{cl6}\gamma (s)\leq \gamma (r).
\ee
It follows from $(\ref{cl5})$, $(\ref{cl6})$ that $\lim_{r\to 0}\gamma (r)=k$. There exists $y_{r_n}$ with $\abs{y_{r_n}}=1$ such that $u(r_ny_{r_n})=\gm_p(r_n)\gamma (r_n)$. Therefore
\be\label{cl14}
\lim_{r_n\to 0}\frac{u(r_ny_{r_n})}{\gm_p(r_n)}=k=\gth.
\ee
Consequently 
\be\label{cl15}
\lim_{r\to 0}\frac{u(ry)}{\gm_p(r)}=\lim_{r\to 0}u_r(y)=\left\{\BA{ll}
k\gm(y)\qquad&\text{if }1<p<N\\
k\qquad&\text{if }p=N.
\EA\right.
\ee
This implies in particular 
\be\label{cl16}
\lim_{x\to 0}\frac{u(x)}{\gm_p(\abs x)}=k.
\ee
Since the convergence of $u_r$ holds in the $C^1_{loc}$-topology, we also deduce that
\be\label{cl17}
\lim_{x\to 0}\abs{x}^{\frac{N-1}{p-1}}\nabla u(x)=\left\{\BA{ll}
\frac{p-N}{p-1}k\frac{x}{\abs x}\qquad&\text{if }1<p<N\\[2mm]
-k\frac{x}{\abs x}\qquad&\text{if }p=N.
\EA\right.
\ee
If we plug these two estimates into the weak formulation of $(\ref{A1})$ we obtain $(\ref{cl2})$.\smallskip
\qeda
\blemma{classlem3}Assume $p$, $q$ are as in \rth{class}. Let $u$ be a positive solution of $(\ref{A1})$ in $\Gw^*$ such that 
\be\label{cl19}
\liminf_{x\to 0}\myfrac{u(x)}{\gm_p(x)}=\infty.
\ee
Then $(\ref{cl3})$ holds.
\es
\Proof If $(\ref{cl4})$ holds, then for any $k>0$ the function $u$ is larger than the radial solution $u_k$ of $(\ref{A1})$ in $B_1^*$ which vanishes on $\prt B_1$ and satisfies $(\ref{cl1})$. When $k\to\infty$ we derive from \rprop{existrad} 
that
\be\label{cl20}
u(x)\geq u_\infty(\abs x)=\gl_{N,p,q}(\abs x^{-\gb_q}-1).
\ee
Next, for any $\ge>0$ we denote by $\tilde u_\ge$ the solution of $(\ref{rad1})$ on $(\ge,1)$ which satisfies 
$\tilde u_\ge(\ge)=\infty$. This solution is expressed from $(\ref{rad13'})$ with a negative $K$, namely
\be\label{cl21}
\tilde u_\ge(r)=b^\frac{1}{q+1-p}\myint{r}{1}s^{\frac{1-N}{p-1}}\left[s^{\frac{q+1-p}{p-1}b}-\ge^{\frac{q+1-p}{p-1}b}\right]^{-\frac{1}{q+1-p}}ds.
\ee
and existence of the blow-up at $r=\ge$ follows from $p>q$. By comparison principle $u\leq\tilde u_\ge+M$ in $B_1\setminus B_\ge$ where $M=\sup\{u(z):\abs z=1\}$. When $\ge\to 0$, formula $(\ref {cl21})$ implies that

\be\label{cl22}\lim_{\ge\to 0}\tilde u_\ge(r)=b^\frac{1}{q+1-p}\myint{r}{1}s^{\frac{1-N}{p-1}}\left[s^{\frac{q+1-p}{p-1}b}\right]^{-\frac{1}{q+1-p}}ds=u_\infty(r).
\ee
Therefore $u_\infty(\abs x)\leq u(x)\leq u_\infty(\abs x)+M$.\qeda\medskip

\noindent{\it Proof of \rth{class}}. By combining \rlemma{classlem1} and \rlemma{classlem3} we have the alternative between (i) and (ii). Assuming now that $\Gw$ is bounded and $u$ and $u'$ are two solutions of $(\ref{A1})$ in $\Gw^*$ continuous in $\overline\Gw\setminus\{0\}$ coinciding on $\prt\Gw$ and satisfying either (i) with the same $k$ or (ii), then, for any $\ge>0$,  $(1+\ge)u+\ge$ is a supersolution which dominates $u'$ in a neighborhood of $0$ and a neighborhood of $\prt\Gw$. Therefore $(1+\ge)u+\ge\geq u'$, which implies $u\leq u'$, and vice versa.\qeda
\medskip

We end this section with a result dealing with global singular solutions.

\bth{class-glob} Let $p-1<q<q_c$ and $1<p\leq N$. If $u$ is a nonnegative solution of $(\ref{A1})$ in $\BBR^{N}_*$, then $u$ is radial and we have the following dichotomy:\smallskip

\noindent (i) either there exists $M\geq 0$ such that $u(x)\equiv M$,

\noindent (ii) either there exist $k>0$, $M\geq 0$ such that $u(x)=u_{k,M}(\abs x)$ defined by $(\ref{rad17})$,

\noindent (ii) or there exists some $M\geq 0$ such that $u(x)=u_{\infty,M}(\abs x)$ defined by $(\ref{rad18})$.
\es
\Proof {\it Step 1: Asymptotic behaviour}. If $u$ is a solution of  $(\ref{A1})$  in an exterior domain $G\supset B_R^c$, it is bounded by \rcor {unif2}. By \rprop{lets0}, it satisfies
\be\label{cl23}
\abs{\nabla u(x)}=(u_r^2+r^{-2}\abs{\nabla'u}^2)^{\frac{1}{2}}(r,\gs)\leq C_{N,p,q}(r-R)^{-\frac{1}{q+1-p}}
\ee
for all $ x=(r,\gs)\in [R,\infty)\ti S^{N-1}$ . Since $q<p$,  
$$\int_{R+1}^\infty\int_{S^{N-1}}\abs {u_r}d\gs dt<\infty,
$$
therefore there exists $\gf\in L^1(S^{N-1})$ such that $u(r,.)\to\gf(.)$ in $L^1(S^{N-1})$ as $r\to\infty$. 
The gradient estimate implies  that the set of functions $\{u(r,.)\}_{r\geq R+1}$ is relatively compact in $C(S^{N-1})$, therefore $u(r,.)\to\gf(.)$ uniformly on $S^{N-1}$ when $r\to\infty$. If $\gs$ and $\gs'$ belong to $S^{N-1}$, there exists 
a smooth path $\gamma:=\{\gamma(t):t\in [0,1]\}$ such that $\gamma (t)\in S^{N-1}$, $\gamma (0)=\gs$, $\gamma (1)=\gs'$. Then
$$u(r,\gs)-u(r,\gs')=\myint{0}{1}\myfrac{d}{dt}u(r,\gamma(t))dt=\myint{0}{1}\langle\nabla'u(r,\gamma(t)),\gamma'(t)\rangle dt,
$$
and finally, using $(\ref{cl23})$, 
\be\label{cl24}
\abs{u(r,\gs)-u(r,\gs')}\leq \norm{\gamma'}_{L^\infty}\abs{\nabla'u(r,\gamma(t))}\leq C_{N,p,q}\norm{\gamma'}_{L^\infty}
r(r-R)^{-\frac{1}{q+1-p}}
\ee
Letting $r\to\infty$, it implies that $\gf$ is a constant, say $M$. As a consequence we have proved that
\be\label{cl25}
\lim_{\abs x\to\infty}u(x)=M.
\ee
Notice that  we did not use the fact that $u$ is a nonnegative solution in order to derive $(\ref{cl24})$. Next we assume the positivity.\smallskip

\noindent{\it Step 2: End of the proof}. If $u$ satisfies $(\ref{cl1})$ for some $k>0$, then for any $\ge>0$, there holds with the notations of \rprop{globrad}
$$(1-\ge)u_{k,M}(\abs x)\leq u(x)\leq (1-\ge)u_{k,M}(\abs x)\qquad\forall x\in\BBR^{N}_*
$$
This implies $u=u_{k,M}$. Similarly, if satisfies $(\ref{cl3})$, we derive $u=u_{\infty,M}$.\qeda

\subsubsection{Negative solutions}
The next result make explicit the behaviour of negative solutions near an isolated singularity.
\bth{classN} Let $p-1<q<q_c$ and $1<p\leq N$. If $u$ is a negative solution of $(\ref{A1})$ in $\Gw\setminus\{0\}$, then there exists $k\leq 0$ such that (\ref{cl1}) and (\ref{cl2}) hold. Furthermore, if $k=0$, $u$ can be extended as a $C^{1,\alpha}$ solution of $(\ref{A1})$ in $\Omega$.
\es
\Proof We can assume $\overline B_1\subset\Gw$. Since $\tilde u:=-u$ satisfies
$$-\Gd_p\tilde u=\abs {\nabla u}^{q}\quad\text{in }\Gw\setminus\{0\}.
$$
It follows from \cite[Th 1.1]{B-V0} that $\abs {\nabla u}^{q}\in L^1_{loc}(\Omega)$ and there exists $k\geq 0$ such that 
\bel{clN1}-\Gd_p\tilde u=\abs {\nabla u}^{q}+k\gd_0\quad\text{in }\CD'(\Gw).
\ee
Furthermore $\abs{\nabla u}^{p-1}\in M^{\frac{N}{N-1}}_{loc}(\Omega)$, where $M^{p}$ denotes the Marcinkiewicz space  (or weak
$L^p$ space). This implies 
$$B:=\abs{\nabla u}^{q+1-p}\in M^{\frac{N(p-1)}{(q+1-p)(N-1)}}_{loc}(\Omega)\subset L^{\frac{N(p-1)}{(q+1-p)(N-1)}-\gs}_{loc}(\Omega)$$ 
for any $\gs>0$. Since 
$q<q_c$, it follows $B\in L^{N+\ge}_{loc}(\Omega)$ for some $\ge>0$. We write the equation under the form
\bel{clN2}
-\Gd_p\tilde u=B\abs{\nabla \tilde u}^{p-1}.
\ee
As a consquence of \cite[Th 1]{Ser} that either there exists $k'>0$ such as 
\bel{clN3}
   \frac{1}{c'}\leq \frac{\tilde u}{\gm_p}\leq c'\quad\text{near }\,0,
\ee
or $u$ has a removable singularity at $0$. If the singularity is removable, then $(\ref{clN1})$ holds with $k=0$. If the singularity is not removable, we set
\bel{clN4}\gamma=\limsup_{x\to 0}\frac{\tilde u(x)}{\gm_p(x)}.
\ee
Then there exists a sequence $\{x_n\}$ converging to $0$ such that 
\bel{clN5}
\gamma=\lim_{n\to \infty}\tilde u(x_n)/\gm_p(x_n)
\ee
 We set $\gd_n=\abs{x_n}$, $\xi_n=x_n/\gd_n$ and
$$\tilde u_{\gd_n}(\xi)=\frac{\tilde u(\gd_n\xi)}{\gm_p(\gd_n)}.
$$
Then 
$$-\Gd_p\tilde u_{\gd_n}-C(\gd_n)\abs{\nabla \tilde u_{\gd_n}}^{q}=0
$$
in $B_{\gd_n^{-1}}\setminus\{0\}$ where
$$C(\gd_n)=\gd_n^{p-q}(\gm(\gd_n))^{q+1-p}.
$$
Since $u_{\gd_n}(\xi)\leq c\gm_p(\xi)$, we derive from \rlemma{classlem0} 
$$\BA {ll}\abs{\nabla u_{\gd_n}(\xi)}\leq c\abs\xi^{-1}\gm_p(\xi)\quad&\text{for } \abs\xi\leq \frac{1}{2\gd_n}\\
\abs{\nabla u_{\gd_n}(\xi)-\nabla u_{\gd_n}(\xi')}\leq c\abs{\xi-\xi'}^\ga\abs\xi^{-1-\ga}\gm_p(\xi)
\quad&\text{for } \abs\xi\leq\abs{\xi'}\leq \frac{1}{2\gd_n}.
\EA$$
Thus, by Ascoli's theorem, the set of functions $\{u_{\gd_n}\}$ is relatively compact in the $C^1_{loc}$-topology of $\BBR^N_*$. Since $C(\gd_n)\to 0$, there exists a subsequence $\{\tilde u_{\gd_{n_k}}\}$ and a nonnegative $p$-harmonic function $\tilde w$ such that $\tilde u_{\gd_{n_k}}\to\tilde v$
as well as its gradient, uniformly on any compact subset of $\BBR^N_*$. 
All the positive $p$-harmonic functions in $\BBR^N_*$ are known (see\cite{KV}: either they are a positive constant, if $N=p$ or have the form  $\gl\gm_p+\gt$ for some $\gl$, $\gt\geq 0$ if $1<p<N$. 

\noindent If $p=N$, we obtain from $(\ref{clN5})$
\bel{clN6}
\lim_{n_k\to\infty}\frac{\tilde u(x_{n_k})}{\gm_p(x_{n_k})}=\gamma=\lim_{n\to \infty}\frac{\tilde u(x_{n})}{\gm_p(x_{n})}
\ee
Thus $\tilde w=\gamma$ and the limit is locally uniform with respect to $\xi$. Therefore for any $\ge>0$, there exists $n_\ge\in\BBN$ such that for $n\geq n_0$, there holds
$$\tilde u(x)\geq (\gamma-\ge)\gm_N(x)\quad\forall x\;\text{ s.t. }\abs x=\gd_n.
$$
By comparison it implies 
$$\tilde u(x)\geq (\gamma-\ge)\gm_N(x)\quad\forall x\;\text{ s.t. }\gd_n\leq\abs x\leq \gd_{n_0}.$$
This holds for any $n\geq n_0$ and any $\ge>0$, therefore,
\bel{clN7}\liminf_{x\to 0}\frac{\tilde u(x)}{\gm_N(x)}\geq \gamma.
\ee
Combining with $(\ref{clN4})$, it implies
\bel{clN8}\lim_{x\to 0}\frac{\tilde u(x)}{\gm_N(x)}= \gamma.
\ee
\noindent If $1<p<N$, estimate $u_{\gd_n}(\xi)\leq C\gm_p(\xi)$ implies $\gt=0$, thus $\tilde w=\gl\gm_p$. Clearly $\gl=\gamma$ because of $(\ref{clN4})$. Similarly as in the case $p=N$, $(\ref{clN7})$ and $(\ref{clN8})$ hold. Since the convergence is in $C^1$, we also get
\bel{clN9}\lim_{x\to 0}\frac{\tilde u_{x_j}(x)}{{\gm_N}_{x_j}(x)}= \gamma.
\ee
From $(\ref{clN1})$ it implies  that there holds
\bel{clN10}
-\Gd_p\tilde u=\abs{\nabla \tilde u}^q+c_{N,p}\gamma\gd_0\quad\text{in }\CD'(\Omega).
\ee
\qeda\medskip

\noindent \Remark In the case $q>q_c$ the description of the isolated singularities is much more difficult, as it is the case if one considers the positive solutions of  
\bel{clN11}
-\Gd_p\tilde u=\tilde u^m\qquad\text{in }\Omega\setminus\{0\}
\ee
for $m>m_c:=\frac{N(p-1)}{N-p}$ (see \cite{SZ} for partial but very deep results). In the case of equation 
\bel{clN12}
-\Gd_p\tilde u=\abs{\nabla \tilde u}^q\quad\text{in }\BBR^N_*
\ee
the main difficulty is to prove that there exists only one positive solution under the form $\tilde u(x)=\tilde u(r,\gs)$, which is the function $\tilde U$. Equivalently it is to prove that the only positive solution of 
\bel{clN13}\BA {ll}
-div\left(\left(\gb_q^2\gw^2+|\nabla'\gw|^2\right)^{\frac{p-2}2}\nabla\gw\right)
-\left(\gb_q^2\gw^2+|\nabla'\gw|^2\right)^{\frac{q}2}\\
\phantom{-div\left(\left(\gb_q^2\gw^2\right)^{\frac{p-2}2}\right)}
-\gb_q(\gb_q(p-1+p-N)(\left(\gb_q^2\gw^2+|\nabla'\gw|^2\right)^{\frac{p-2}2}\gw
=0\quad\text{in }S^{N-1}
\EA\ee
is the constant $\tilde \gl_{N,p,q}$.
\mysection{Quasilinear equations on Riemannian manifolds}
\subsection{Gradient geometric estimates}

In this section we assume that  $(M^{^{_N}}\!\!,g)$ is a N-dimensional Riemannian manifold, $TM$ its tangent bundle, $\nabla u$ is the covariant gradient, $\langle.,.\rangle$ the scalar product expressed in the metric $g:=(g_{ij})$, $Ricc_g$ the Ricci tensor and $Sec_g$ the sectional curvature. Formula $(\ref{W1})$ is a particular case of the B\"ochner-Weitzentb\"ock formula which is the following: if $u\in C^3(M)$ there holds
\bel{W2}
\myfrac{1}{2}\Gd_g\abs{\nabla u}^2=\abs{D^2u}^2+\langle\nabla \Gd_g u,\nabla u\rangle+Ricc_g(\nabla u,\nabla u), 
\ee
where $D^2u$ is the Hessian, $\Gd_g={\rm div}_g(\nabla u)$ is the Laplace-Beltrami operator on $(M^{^{_N}}\!\!,g)$ and ${\rm div}_g$ is the divergence operator acting on $C^1(M,TM)$. For $p>1$, we also denote by $\Gd_{g,p}$ the p-Laplacian operator on $M$ defined by
\bel{p-lapl}
\Gd_{g,p}u={\rm div}_g(\abs{\nabla u}^{p-2}\nabla u),
\ee
with the convention $\Gd_{2,g}=\Gd_g$. A natural geometric assumption is that the Ricci curvature is bounded from below and more precisely
\bel{W3}Ricc_g(x)(\xi,\xi)\geq -(N-1)B^2\abs\xi^2\qquad\forall \xi\in T_xM
\ee
for some $B\geq 0$. If $u\in C^3(M)$ is a solution  of 
\bel{S1''}
-\Gd_{p,g}u+\abs{\nabla u}^{q}=0\qquad\text{in }M,
\ee
then $(\ref{grad4})$ is replaced by
\bel{W4}\BA {l}
\Gd_g z+(p-2)\myfrac{\langle D^2z(\nabla u),\nabla u\rangle}{z}\geq\myfrac{2a^2}{N}z^{q+2-p}-\myfrac{1}{Na^2}\myfrac{\langle\nabla z,\nabla u\rangle^2}{z^2}-\myfrac{(p-2)}{2}\myfrac{\abs{\nabla z}^2}{z}
\\[3mm]\phantom{\Gd z----}
+(p-2)\myfrac{\langle \nabla z,\nabla u\rangle^2}{z^2}+
(q+2-p)z^{\frac{q-p}{2}}\langle\nabla z,\nabla u\rangle-(N-1)B^2z.
\EA\ee
and $\CL$ in $(\ref{grad7})$ by 
\bel{W5}
\CL^*(z):=\CA (z)+Cz^{q+2-p}- D\myfrac{\abs{\nabla z}^2}{z}-(N-1)B^2z\leq 0\qquad\text {in }\Gw.
\ee

We recall that the convexity radius $r_M(a)$ of some $a\in M$ is the supremum of all the $r>0$ such that the ball $B_r(a)$ is convex. 
Note that, in order to obtain estimates on the gradient of solution, when $p\neq 2$ an extra assumption besides $(\ref{W3})$ is needed; it concerns the sectional curvature. 
\blemma{compa} Assume $q>p-1\geq 0$ and let $a\in M$, $R>0$ and $B\geq 0$ such that $Ricc_g\geq -(N-1)B^2$ in $B_R(a)$. Assume also $Sec_g\geq -\tilde B^2$ in $B_R(a)$ for some $\tilde B\geq 0$ if $p>2$, or 
$r_M(a)\geq R$ if $1<p<2$. Then there exists $c=c(N,p,q)>0$ such that the function 
\bel{Ge0}w(x)=\gl\left(R^2-r^2(x)\right)^{-\frac{2}{q+1-p}}+\gm,
\ee
where $r=r(x)=d(x,a)$, satisfies
\bel{Ge1}
\CL^*(w)\geq 0 \qquad\text{in }B_R(a),
\ee
provided that
\bel{Ge2}\BA{l}
\gl=c\max\left\{(R^4B^2)^{\frac{1}{q+1-p}},((1+(B+(p-2)_+\tilde B)R)R^2)^{\frac{1}{q+1-p}}\right\}
\EA\ee
and
\bel{Ge2'}
\gm\geq ((N-1)B^2)^{\frac{1}{q+1-p}}.
\ee
\es
\Proof Let $w$ as in $(\ref{Ge0}))$. We will show that by choosing $\gl$ and $\gm$ as in $(\ref{Ge2}))$ and $(\ref{Ge2'})$
respectively, then $(\ref{Ge1}))$ holds.
We recall that
\bel{Ge3}
\Gd_gw=w''+w'\Gd_gr,
\ee
and (see \cite[Lemma 1]{RRV} )
$$\Gd_g r\leq (N-1)B\coth ( Br)\leq \myfrac{N-1}{r}\left(1+Br\right).$$ 
Then
\bel{Ge4}\BA {ll}
\Gd_gw\leq \myfrac{4}{q+1-p}(R^2-r^2)^{-\frac{2(q+2-p)}{q+1-p}}\\[4mm]
\phantom{--------}
\displaystyle \ti\left(\frac{2r^2(q+3-p)}{q+1-p}+(R^2-r^2)(1+(N-1)\left(1+Br\right)\right).
\EA\ee
Moverover \cite[Chapt 2, Theorem A]{GW}
\bel{Ge5}
D^2w=w''dr\otimes dr+w'D^2r.
\ee
If $r_M(a)\geq r(x)$, then the ball $B_{r(x)}(a)$ is convex. This implies that $r$ is convex and therefore $D^2r\geq 0$ (see \cite[IV-5]{Sa}).
Furthermore, if $Sec_g(x)\geq -\tilde B^2$, then from \cite[IV-Lemma 2.9]{Sa},
\bel{Ge6}
D^2r\leq \tilde B\coth (\tilde Br)g_{ij}\leq \myfrac{\tilde B}{r}\left(1+\tilde Br\right)g_{ij},
\ee
therefore
\bel{Ge7}\BA {ll}
0\leq \myfrac{\langle D^2w(\nabla u),\nabla u\rangle}{\abs{\nabla u}^2}\\[3mm]
\phantom{0}\leq 
\myfrac{4}{q+1-p}(R^2-r^2)^{-\frac{2(q+2-p)}{q+1-p}}\left(\myfrac{2r^2(q+3-p)}{q+1-p}+(R^2-r^2)(2+\tilde Br)\right).
\EA\ee
We obtain
\bel{Ge8}
\BA {l}\CA(w)= -\Gd w-(p-2)\myfrac{\langle D^2w(\nabla u),\nabla u\rangle}{\abs{\nabla u}^2}\\[2mm]
\phantom{\CA(w)}
\geq -k\gl(R^2-r^2)^{-\frac{2(q+2-p)}{q+1-p}}(R^2+(p-2)_+(R^2-r^2)\tilde Br\coth \tilde Br)\\[2mm]
\phantom{\CA(w)}
\geq  -k\gl(R^2-r^2)^{-\frac{2(q+2-p)}{q+1-p}}(R^2+(R^2-r^2)B_pr)
\EA\ee
for some $k=k(N,p,q)$, where $B_p=B+(p-2)_+\tilde B$. Since
$$w^{q+2-p}\geq \gl^{q+2-p}\left(R^2-r^2\right)^{-\frac{2(q+1-p)}{q+1-p}}+\gm^{q+2-p},
$$
we have
\bel{Ge9}\BA{l}
\CL^*(w)\\\phantom{--}
\geq \gl(R^2-r^2)^{-\frac{2(q+2-p)}{q+1-p}}\left(-k(R^2+(R^2-r^2)B_pr)-\myfrac{16D}{(q+1-p)^2}r^2
+C\gl^{q+1-p}\right)\\[4mm]\phantom{\CL^*(w)}
+\gm^{q+2-p}-(N-1)B^2\gl\left(R^2-r^2\right)^{-\frac{2}{q+1-p}}-(N-1)B^2\gm.
\EA\ee
Take 
$\gm\geq ((N-1)B^2)^{\frac{1}{q+1-p}}
$ as in $(\ref{Ge2'})$.
Next we choose $\gl$ in order to have, uniformly for $0\leq r<R$,
$$
2^{-1}C\gl^{q+1-p}\geq k(R^2+(R^2-r^2)B_pr)+D\myfrac{16}{(q+1-p)^2}r^2,
$$
so that
$$\BA {ll}
\CL^*(w)\geq 2^{-1}C\gl(R^2-r^2)^{-\frac{2(q+2-p)}{q+1-p}}\gl^{q+1-p}-(N-1)B^2\gl(R^2-r^2)^{-\frac{2}{q+1-p}}\\
[2mm]\phantom{\CL^*(w)}
=\gl(R^2-r^2)^{-\frac{2(q+2-p)}{q+1-p}}\left(2^{-1}C\gl^{q+1-p}-(N-1)B^2(R^2-r^2)^2\right)
\EA$$
uniformly for $0\leq r<R$. Then we enlarge $\gl$ if necessary to have
$$2^{-1}C\gl^{q+2-p}(R^2-r^2)^{-\frac{2(q+2-p)}{q+1-p}}\geq (N-1)B^2\gl\left(R^2-r^2\right)^{-\frac{2}{q+1-p}},
$$
also uniformly for $0\leq r<R$. Hence we see that there exists $c=c(N,p,q)$ such that, if we choose
\bel{Ge11}\BA{l}
\gl=c\max\left\{(R^4B^2)^{\frac{1}{q+1-p}},((1+B_pR)R^2)^{\frac{1}{q+1-p}}\right\}
\EA\ee
then $(\ref{Ge1})$ holds.\qeda

\bprop{up} Assume $q>p-1>0$. Let $\Gw$ be an open subset of $M$ such that $Ricc_g\geq (1-N)B^2$ in $\Gw$. Assume also $Sec_g\geq -\tilde B^2$ in  $\Gw$ if $p>2$, or  $r_M(x)\geq \dist (x,\prt\Gw)$ for any $x\in M$ if $1<p<2$. Then any solution $u$ of $(\ref{S1''})$ in $\Gw$ satisfies
\bel{est1}
\abs{\nabla u(x)}^2\leq c_{N,p,q}\max\left\{B^{\frac{2}{q+1-p}},(1+B_pd(x,\prt\Gw)))^{\frac{1}{q+1-p}}(d(x,\prt\Gw))^{-\frac{2}{q+1-p}}\right\}\quad\forall x\in \Gw,
\ee 
where $B_p=B+(p-2)_+\tilde B$.
\es
\Proof Assume $a\in \Gw$ and $R<d(a,\prt\Gw)$. Let $w$ be as in \rlemma{compa}, then
\bel{Ge12}
\CA (z-w)+C\left(z^{q+2-p}-w^{q+2-p}\right)-(N-1)B^2(z-w) -D\left(\myfrac{\abs{\nabla z}^2}{z}-\myfrac{\abs{\nabla w}^2}{w}\right)\leq 0
\ee
in $B_R(a)$. Let $G$ be a connected component of the set $\{x\in B_R(a):z(x)-w(x)>0\}$. Then, if $C(q+2-p)(w(a))^{q+1-p}>(N-1)B^2$, by the mean value theorem and the fact that $w(a)$ is the minimum of $w$, there holds that 
\bel{Ge12-1}
C\left(z^{q+2-p}-w^{q+2-p}\right)-(N-1)B^2(z-w)>0 \qquad\text{in }G.
\ee
Since $w(a)\geq\gm\geq ((N-1)B^2)^{\frac{1}{q+1-p}}$ and $q+2-p>1$, this condition is fulfilled by choosing the right $\gm$ as in $(\ref{Ge2'})$. We conclude as in the proof of \rprop{lets0} that $G=\emptyset$. Therefore $z\leq w$ in $B_R(a)$. In particular,  
\bel{Ge13}
z(a)\leq c_{N,p,q}\max\{B^{\frac{2}{q+1-p}},(1+B_pR)^{\frac{1}{q+1-p}}R^{-\frac{2}{q+1-p}}\}
\ee
where $c_{N,p,q}>0$. Then $(\ref{est1})$ follows.\qeda\medskip

\noindent \Remark Since $Ricc_g(x)(\xi,\xi)=(N-1)\sum_VSec_g(x)(V)$, where $V$ denotes the set of two planes in $T_xM$ which contain $\xi$, there holds
 $$Sec_g\geq -\tilde B^2\Longrightarrow Ricc_g\geq (1-N)\tilde B^2.$$
  However, in the previous estimate, the long range estimate on $\nabla u$ depends only on the Ricci curvature.\medskip
    \subsection{Growth of solutions and Liouville type results}
\bcor{liouville2} Assume $(M^{^{_N}}\!\!,g)$ is a complete noncompact N-dimensional Riemannian manifold such that $Ricc_g\geq (1-N)B^2$  and let $q>p-1>0$. Assume also if $r_M(x)=\infty$ if $1<p<2$ or that the sectional curvature $Sec_g$  satisfies for some $a\in M$
\bel{scal}
\lim_{\dist (a,x)\to\infty}\myfrac{\abs{Sec_g(x)}}{\dist (a,x)}=0,
\ee
if $p>2$. Then any solution $u$ of $(\ref{S1''})$ satisfies
\bel{est3-}
\abs{\nabla u(x)}^2\leq c_{N,p,q}B^{\frac{2}{q+1-p}}\qquad\forall x\in M.
\ee
In particular, $u$ is constant if $Ricc_g\geq 0$, while in the general case $u$ has at most a linear growth with respect to the distance function.  
\es

\noindent {\it Application} An example of a complete manifold with constant negative Ricci curvature is the standard hyperbolic space $(\BBH^N,g_{_0})$ for which $Ricc_{g_{_0}}=-(N-1)g_0$. Another application deals with positive p-harmonic functions (for related results  with $p=2$  see \cite{Yau}, \cite{CoMin}).
\bcor{liouville3} Assume $(M^{^{_N}}\!\!,g)$ is as in \rcor{liouville2}. Let $p>1$ and assume that $(\ref{scal})$ holds if $p> 2$ or $Sec_g\leq 0$ if $1<p<2$. If  $v$ is a positive p-harmonic function, then \smallskip

\noindent (i)  if $Ricc_g\geq 0$, $v$  is constant. \smallskip

\noindent (ii) if $\inf\{Ricc_g(x):x\in M\}=(1-N)B^2<0$, $v$ satisfies
\bel{est3-'}
 v(a) e^{-c_{N,p}B\dist(x,a)}\leq v(x)\leq v(a) e^{c_{N,p}B\dist(x,a)}\qquad\forall x\in M.
\ee
\es
\Proof We take $q=p$ and assume that $v$ is p-harmonic and positive. If we write $v=e^{-\frac{u}{p-1}}$ , then $u$ satisfies 
$$-\Gd_{g,p} u+\abs{\nabla u}^p=0.$$
If $Ricc_g(x)\geq 0$, $u$, and therefore $v$ is constant by \rcor{liouville2}. If $\inf\{Ricc_g(x):x\in M\}=(1-N)B^2<0$
we apply $(\ref{est3-})$ to $\nabla u$. If $\gamma$ is a minimizing geodesic from $a$ to $x$, then $\abs{\gamma'(t)}=1$
and 
$$u(x)-u(a)=\int_{0}^{d(x,a)}\myfrac{d}{dt}u\circ\gamma (t) dt=\int_{0}^{d(x,a)}\langle\nabla u\circ\gamma (t),\gamma' (t)\rangle dt.
$$
Since
$$\abs{\langle\nabla u\circ\gamma (t),\gamma' (t)}\leq \abs{ \nabla u\circ\gamma (t)}\leq c_{N,p}\gk,
$$
we obtain
\bel{est3}
u(a)-c_{N,p}B\dist(x,a)\leq u(x)\leq u(a)+c_{N,p}B\dist(x,a)\qquad\forall x\in M,
\ee
Then $(\ref{est3-'})$ follows since $u=(1-p)\ln v$. Notice that (i) follows from (ii) and that in the case $1<p<2$ the assumption (i) implies that $Ricc_g= 0$.\qeda



\begin{thebibliography}{99}



\bibitem{BP} Baras P., Pierre M., {\textit Singularit\'{e}s \'{e}liminables pour des \'{e}quations semi-lin\'{e}aires}, Ann. Inst. Fourier {\bf  34}, 185-206 (1984).


\bibitem{B-V0} Bidaut-V\'{e}ron M.F., \textit{Local and global behavior of solutions of quasilinear equations of Emden-Fowler type}, Arch. Rat. Mech. Anal. {\bf 107}, 293-324 (1989).

\bibitem{B-V} Bidaut-V\'{e}ron M.F., \textit{Removable singularities and
existence for a quasilinear equation}, Adv. Nonlinear Studies {\bf 3}, 25-63 (2003).


\bibitem{B-VH} Bidaut-V\'{e}ron M.F., Hamid A. \textit{On the connection between two quasilinear ellpitic problems with lower order terms of order 0 or 1}, Comm. Cont. Math. {\bf 12}, 727-788 (2010).

\bibitem{BV-H-V} Bidaut-V\'{e}ron M.F., Nguyen Quoc H., V\'eron L. \textit{Quasilinear Lane-Emden equations with absorption 
and measure data}, J. Math. Pures Appl. http://dx.doi.org/10.1016.matpur.2013.11.011, to appear.

\bibitem{BoGa} Boccardo L., Gallouet T., \textit{Nonlinear elliptic and
parabolic equations involving measure data,} J Funct. Anal. {\bf 87},
149-169 (1989).

\bibitem{BoGa1} Boccardo L., Gallouet T., \textit{Nonlinear elliptic  equations with right-hand side measures}, Comm. Part. Diff. Eq. {\bf 17}, 641-655 (1992).






\bibitem{CY} S. Y. Cheng and S.-T. Yau: {\em Differential equations on Riemannian manifolds and their geometric applications}, { Comm. Pure Appl. Math. 28}, 333-354 (1975).


\bibitem{CoMin} Colding T. H., Minicozzi II W. P., \textit{Harmonic functions on manifolds}, 
Annals of Mathematics {\bf 146}, 725-747 (1997).


\bibitem{DMOP} Dal Maso G., Murat F., Orsina L., Prignet A., \textit{%
Renormalized solutions of elliptic equations with general measure data,}
Ann. Scuola Norm. Sup. Pisa {\bf 28}, 741-808 (1999).


\bibitem{FaSe} Farina A., Serrin J., \textit{Entire solutions of completely coercive quasilinear elliptic equations, II }, J. Diff. Eq. {\bf 250}, 4367-4408 (2011).

\bibitem{FV} Friedman A., V\'eron L. \textit{Singular Solutions of Some
Quasilinear Elliptic Equations,}, Arch. Rat. Mec. Anal. {\bf 96},
258-287 (1986).




\bibitem{GW} Greene R.E., Wu H., \textit {Function Theory on Manifolds
Which Possess a Pole}, Lecture Notes in Mathematics {\bf 699}, Sringer-Verlag Berlin Heidelberg (1979).

\bibitem{HMV} Hansson K., Maz'ya V.G., Verbitsky I. E., \textit{ Criteria of solvability for multidimensional Riccati equations}, Ark. Math. {\bf 37}, 87-120 (1999).

\bibitem{HKM} Heinonen J., Kilpelainen T., Martio O., \textit{ Nonlinear potential theory of degenerate elliptic equations, } Oxford Science Publ. (1993).


\bibitem{Iv} Iwaniec T., \textit {Projection onto gradient fields and $L^p$ estimates for degenerate elliptic operators, }Studia Math. {\bf 75}, 293-312 (1983).

\bibitem{KV} Kichenassamy S., V\'eron L. \textit{Singular solutions of the p-Laplace equation,} Math. Ann. {\bf 275}, 599-615 (1986).

\bibitem{KKTK} Kilpelainen T., Kuusi T., Tuhuola-Kujanpaa A.  \textit{Superharmonic functions are locally renormalized}, Ann. I. H. P. Anal. Non Lin\'eaire {\bf 28}, 775-795 (2011).

\bibitem{KL} Kilpelainen T., Li, G., \textit{Estimates for $p$-Poisson equations}, Diff. Int. Equ. {\bf 13}, 781-800 (2000).

\bibitem{KZ} Kilpelainen T., Zou J., \textit{ A local estimate for nonlinear equations with discontinuous coefficients}, Comm. Part. Diff. Eq. {\bf 24}, 2043-2068 (1999).

\bibitem{KZ1} Kilpelainen T., Zou J., \textit{ A boundary estimate for  nonlinear equations with discontinuous coefficients}, Diff. int. Eq. {\bf 14}, 475-492 (2001).


\bibitem{KN} Kotschwar B. L., Ni L., \textit{  Local gradient estimates of $p$-harmonic functions, $1/H$-flow, and an entropy formula},   Ann. Sci. \'Ecole. Norm. Sup.  {\bf 42}, 1-36 (2009).

\bibitem{LaLi} Lasry J.M., Lions P. L., \textit{ Nonlinear elliptic equations with singular boundary conditions and stochastic control with state constraints. I. The model problem}, Math. Ann. {\bf 283}, 583-630(1989) .

\bibitem{Lio} Lions P. L., \textit{Quelques remarques sur les probl\`emes elliptiques quasilin\'eaires du second ordre,} J. Anal.Math.{\bf  45}, 234-254  (1985).


\bibitem{Mae} Maeda F., \textit{Renormalized solutions of Dirichlet problems
for quasilinear elliptic equations with general measure data}, Hiroshima
Math. J., {\bf 38}, 51-93 (2008).

\bibitem{MuPor} Murat F., Porretta A., \textit{Stability properties,
existence, and nonexistence of renormalized solutions for elliptic equations
with measure data,} Comm. Partial Differential Equations {\bf 27}, 2267-2310 (2002).

\bibitem{Phuc} Phuc N. C. \textit {Riccati type equations with super-critical exponents.},  Comm. Part. Diff. Eq. {\bf 35}, 1958-1981 (2010).


\bibitem{Phuc0} Phuc N. C. \textit {Nonlinear Muckenhoupt-Wheeden type bounds on Reifenberg flat domains, with applications to quasilinear Riccati type equations.},  Adv. Math. {\bf 250}, 387-419 (2014).

\bibitem{Phuc1} Phuc N. C. \textit {Global integral gradient bounds for quasilinear equations below or near the critical exponent}, Ark. Math., doi: 10.1007/s11512-012-0177-5, to appear.

\bibitem{PhV1} Phuc N. C. , Verbitsky I.E.  \textit {Quasilinear and Hessian equations of Lane-Emden type,} Ann. Math. {\bf 168}, 859-914 (2008).

\bibitem{PhV2} Phuc N. C. , Verbitsky I.E.  \textit { Singular quasilinear and Hessian equations and inequalities}, J. Funct. Anal. {\bf 259}, 1875-1906 (2009).

\bibitem{NpV} Nguyen Phuoc T., V\'eron L. {\textit Boundary singularities of solutions to elliptic viscous
Hamilton-Jacobi equations}, J. Funct. An. {\bf 263}, 1487-1538  (2012).


\bibitem{PPP} Petitta F., Ponce A., Porretta A. \textit { Diffuse measures and nonlinear parabolic equations}
J. Evol. Equ. {\bf 11}, 861?905 (2011).






\bibitem{RRV} Ratto A., Rigoli M., V\'eron L.  \textit{Conformal immersion of complete Riemannian manifolds and extensions of the Schwarz lemma}, Duke Math. J. {\bf 74}, 223-236 (1994).

\bibitem{Sa} Sakai T., \textit {Riemannian Geometry}, Transl. Math. Mono. {\bf 149}, Amer. Math. Soc. (1997). 


\bibitem{Ser0} Serrin J.\textit{ Local behaviour of solutions of 
quasilinear equations}, Acta Math. {\bf 111}, 247-302 (1964).


\bibitem{Ser} Serrin J. \textit {Isolated singularities of solutions of 
quasilinear equations}, Acta Math. {\bf 113}, 219-240 (1965).

\bibitem{SZ} Serrin J., Zou H. \textit{ Cauchy-Liouville and universal boundedness theorems for quasilinear elliptic equations}
Acta Math. {\bf 189}, 79-142 (2002).

\bibitem{To} Tolksdorff P. \textit{ Regularity for a more general class of quasilinear elliptic equations}, J. Diff. Equ. {\bf 51}, 
126-150 (1984).

\bibitem{Tru} Trudinger N., \textit{On Harnack type inequalities and their applications to quasilinear elliptic equations},  Comm. Pure Appl. Math. {\bf 20}, 721-747 (1967).


\bibitem{Yau} Yau S. T., \textit{Harmonic functions on complete Riemannian manifolds}, Comm. Pure and Appl. Math. {\bf 28}, 201-228 (1975).
\end{thebibliography}
\end{document}